\def\E{\end{document}}
\documentclass[11pt]{article}
\usepackage{}
\usepackage{amssymb}
\usepackage{amssymb,amsmath,color}

\begin{document}
\title{\bf Global small solutions to heat conductive compressible nematic liquid crystal system: smallness on a
scaling invariant quantity}

\author{Jinkai Li \thanks {
South China Research Center for Applied Mathematics and Interdisciplinary Studies,
School of Mathematical Sciences, South China Normal University, Guangzhou 510631, China.
 E-mail: jklimath@m.scnu.edu.cn; jklimath@gmail.com} \quad and  \quad Qiang Tao\thanks {Corresponding author. School of Mathematics and Statistics,
 Shenzhen University, Shenzhen 518060, China and Shenzhen Key Laboratory of Advanced Machine Learning and Applications, Shenzhen University, Shenzhen 518060, China. E-mail: taoq060@126.com},
%\and
}

\date{}
\maketitle
\newtheorem{theorem}{Theorem}[section]
\newtheorem{definition}{Definition}[section]
\newtheorem{lemma}{Lemma}[section]
\newtheorem{proposition}{Proposition}[section]
\newtheorem{corollary}{Corollary}[section]
\newtheorem{remark}{Remark}[section]
\renewcommand{\theequation}{\thesection.\arabic{equation}}
\catcode`@=11 \@addtoreset{equation}{section} \catcode`@=12

{\bf Abstract.}\ \ {\small
In this paper, we consider the Cauchy problem to the three dimensional heat conducting
compressible nematic liquid crystal system in the presence of vacuum and with vacuum far fields. Global
well-posedness of strong solutions is established under the condition that the scaling invariant quantity
$ (\|\rho_0\|_\infty+1)\big[\|\rho_0\|_3+(\|\rho_0\|_\infty+1)^2(\|\sqrt{\rho_0}u_0\|_2^2+ \|\nabla d_0\|_2^2)\big]
\big[\|\nabla u_0\|_2^2+(\|\rho_0\|_\infty+1)(\|\sqrt{\rho_0}E_0\|_2^2 + \|\nabla^2 d_0\|_2^2)\big]$
is sufficiently small with the smallness depending only on the parameters appeared in the system.}


%{\bf subject classifications.}\ \  {35B40, 35Q35, 76N10.}

\section{Introduction}
Liquid crystals are intermediate phases between solids and
fluids. The continuum theory of
liquid crystals was established by Ericksen \cite{E1962} and Leslie \cite{L1968} during the period of 1958 through
1968. The present paper concerns a simplified version of the general Ericksen-Leslie system,
which roughly speaking is a coupled system of the compressible Navier-Stokes equations and the
harmonic heat flow (see \cite{L1989, LL1995}). The equations of the heat conducting compressible nematic
liquid system read as:
\begin{align}
&\rho_t +{\rm div}(\rho u)=0, \label{1.1}
\\
&\rho(u_t+ u\cdot\nabla u)+\nabla P=\mu \Delta  u+
(\mu+\lambda)\nabla {\rm div}u- \Delta d\cdot\nabla d,
\label{1.2}\\
&c_v\rho(\theta_t+u\cdot\nabla\theta)+P{\rm div}u-\kappa\Delta\theta=\mathcal{Q}(\nabla u)+|\Delta d+|\nabla d|^2d|^2,
\label{1.3}\\
&d_t+u\cdot\nabla d=\Delta d+|\nabla d|^2d, \label{1.4}
\end{align}
%%%%
%%
where $\rho \in \mathbb{R}_+$ is the density of the fluid, $u \in \mathbb{R}^3$ is the velocity field and $d \in \mathbb{S}^2 $ represents
macroscopic average of the nematic liquid crystal orientation field, with $\mathbb{S}^2 $ denoting the unit
spherical surface in $\mathbb{R}^3$.
Here, $P= R\rho\theta$ is the pressure with $R$ being a positive constant,
$\lambda$ and $\mu$ are constant viscosity coefficients satisfying the physical conditions $\mu>0$ and $2\mu+3\lambda\geq 0,$
heat capacity $c_v=\frac{R}{\gamma-1}$ with $\gamma>1$ being the adiabatic constant, $\kappa>0$ is the heat conductive
constant coefficient, and
$$\mathcal{Q}(\nabla u)= \frac{\mu}{2}|\nabla u + (\nabla u)^t|^2+ \lambda({\rm div}u)^2.$$
The additional assumption that $2\mu>\lambda$ will also be used in this paper.

Mathematical analysis of the nematic liquid crystals have attracted a lot of attentions for
several decades. For the incompressible case, Lin \cite{L1989} first introduced and studied a simplified Ericksen-Leslie system modeling the incompressible liquid crystal flows. From then on, the solvability
and stability of the incompressible liquid crystal flows have been substantially developed. The global existence and partial regularity of weak solutions to the Ginzburg-Landau type approximation system were obtained by Lin and Liu \cite{LL1995, LL1996}, while the global existence of weak solutions to the original system was proved by Lin et al.\,\cite{LINLINWANG}, Hong \cite{HONGWEAK}, and Hong and Xin \cite{HONGXINWEAK} for two-dimensional case, and by Lin and Wang \cite{LINWANGWEAK3D} for
three-dimensional case under some geometric assumptions on the initial director field $d_0$. Uniqueness of weak solutions
in two dimensions was proved by Lin and Wang \cite{LINWANGUNIQUENESS}, Li et al.\,\cite{LITITIXIN}, and Wang et al.\,\cite{WANGWANGZHANGUNIQUENESS}, while the nonuniqueness of weak solutions and finite time blow up of classic solutions in three dimensions was, respectively, 
addressed by Gong et al.\,\cite{NONUNIQUENESS} and Huang et al.\,\cite{BLOWUP}.
Feireisl et al.\,\cite{FFRS2012, FRS2011} considered a non-isothermal Ginzburg-Landau
model of nematic liquid crystals and investigated the global existence of weak solutions, while the global existence of weak solutions to the corresponding non-isothermal Ericksen-Leslie system in two dimensions was established by Li and Xin \cite{LIXINNONISO}. One can refer to De Anna and Liu \cite{DeAnnaLiu} for the derivation of the general non-isothermal
Ericksen-Leslie system.
Huang and Wang \cite{HW2012} established a blow-up criterion
for the short time classical solutions to incompressible liquid crystal flows in dimensions two and three.
Hong et al. \cite{HLZ2014} established the local well-posedness and blow-up criteria of strong solutions
to the liquid crystal system with general Oseen-Frank free energy density.

Concerning the compressible case, the model of the
liquid crystals becomes more complicated since the density variation affects the mechanical
behaviour of the fluid. Global well-posedness of the isentropic compressible nematic liquid
crystals in one dimension was proved by Ding et al. \cite{DWW2011, DLWW2012}, while the global existence of weak
solutions in multi-dimensions was proved by Jiang et al. \cite{JJW2013, JJW2014} under a smallness condition on
the third component of initial orientation field. The local existence of unique strong solution to
the initial value or initial-boundary value problem was proved in Huang et al. \cite{HWW2012, HWW2012-2},
where a series of blow-up criterion of strong solutions were established as well. Li et al. \cite{LXZ2018} obtained the global classical solutions to the Cauchy problem with small initial energy but possibly large oscillations and the initial density may allowed vacuum. The long-time behavior of classical solution was considered in \cite{GTY2016}.
By virtue of the Fourier splitting method, the authors built optimal temporal decay rate of the global solution.
For more results on simplified isothermal Ericksen-Leslie system, the readers can refer to \cite{HW2013, GTY2015, LLW2015, W2016}  and
references therein.

Inspired by the introduction of non-isothermal models of incompressible nematic liquid crystals by Feireisl et al. in \cite{FFRS2012, FRS2011},
the compressible non-isothermal nematic liquid crystal flows is now attracting increasing research attention.
Fan et al. \cite{FLN2018} first investigated the local existence of unique strong solution to the initial boundary
value problem. Guo et al. \cite{GXX2017} obtained the global existence of smooth solutions for the Cauchy problem provided that the initial datum is close to a steady state and
gave the algebraic decay rate of the global solution. A blow up criterion was established in \cite{ZX2019} for the strong solutions to the two-dimensional non-isothermal flows in a bounded domain under a geometric condition introduced by Lei et al. in \cite{LLZ2014}.
Recently, Liu and Zhong \cite{LZ2021}
proved the global well-posedness of strong solutions with the initial data can have compact support provided that the quantity $\|\rho_0\|_{L^\infty}+ \|\nabla d_0\|_{L^3}$ is suitably small with the
smallness depending not only on the parameters involved in the system, but also on some high order norm of the initial data.

The purpose of this paper is to investigate the global well-posedness of strong solutions to the Cauchy problem of (\ref{1.1})-(\ref{1.4})
along with the initial condition:
\begin{align}
\label{1.5}
(\rho, u, \theta, d)|_{t=0}=(\rho_0, u_0, \theta_0, d_0).
\end{align}
The initial data allows far field vacuum and the smallness assumption, which depends only on the parameters involved in the system, is imposed exclusively
on some quantities that are scaling invariant with respect to the following scaling
transform:
\begin{align}
\label{1.6}
(\rho_{0\lambda}(x), u_{0\lambda}(x), \theta_{0\lambda}(x), d_{0\lambda}(x)) =(\rho_0(\lambda x), \lambda u_0(\lambda x), \lambda^2\theta_0(\lambda x), d_0(\lambda x)), \text{~~for any~~} \lambda\neq0.
\end{align}
This scaling transform on the initial data is motivated from the natural scaling invariant property of
the compressible nematic liquid crystal flow (\ref{1.1})-(\ref{1.4}):
\begin{align*}
\rho_\lambda(x,t)=\rho(\lambda x, \lambda^2t), ~u_\lambda(x,t)=\lambda u(\lambda x, \lambda^2t),
~\theta_\lambda(x,t)=\lambda^2\theta(\lambda x, \lambda^2t),~d_\lambda(x,t)=d(\lambda x, \lambda^2t).
\end{align*}
More precisely, if $(\rho, u, \theta, d)$ is a solution with the initial data $(\rho_0, u_0, \theta_0, d_0)$, then,
by straightforward calculations, one can find that $(\rho_\lambda, u_\lambda, \theta_\lambda, d_\lambda)$
is also a solution with the transformed initial data $(\rho_{0\lambda}, u_{0\lambda}, \theta_{0\lambda}, d_{0\lambda})$
for any nonzero $\lambda$.

As already explained in \cite{L2020}, imposing smallness assumptions on the scaling invariant quantities is necessary for obtaining globally well-posed to system (\ref{1.1})-(\ref{1.5}). In fact,
if assuming that $\mathcal{M}$ is a functional such that
$$\mathcal{M}(\rho_{0\lambda}, u_{0\lambda}, \theta_{0\lambda}, d_{0\lambda})=\lambda^\ell(\rho_0, u_0, \theta_0, d_0), ~~\text{for any~} \lambda\neq 0 ~~\text{and some constant~} \ell\neq0,$$
and
that the global well-posedness holds for any initial data $(\rho_0, u_0, \theta_0, d_0)$ satisfying $$\mathcal{M}(\rho_0, u_0, \theta_0, d_0)\leq \varepsilon_0,$$
then, by suitably choosing the scaling parameter $\lambda$, one can show that the global well-posedness holds for arbitrary large initial data. However, such
global well-posedness for arbitrary large initial data is far from what we have already known.

Throughout this paper, the following notations are needed.
For $ 1\leq p\leq\infty $, denote $L^p=L^p(\mathbb{R}^3)$ as the standard $L^p$ Lebesgue spaces with the norm $\|\cdot\|_{p}$.
 For $ 1\leq p\leq\infty $ and positive integer $k$, denote by $W^{k,p}=W^{k,p}(\mathbb{R}^3)$ the standard Sobolev spaces, whose norm is denoted as $\|\cdot\|_{W^{k,p}}$ or $\|\cdot\|_{H^{k}}$ with $H^k=W^{k,2}$.
To simplify the expressions, the norm $\sum^K_{i=1}\|f_i\|_X$ or $\left(\sum^K_{i=1}\|f_i\|_X^2\right)^{\frac{1}{2}}$ are sometimes denoted by $\|(f_1, f_2, \ldots, f_K)\|_X$.
For $ 1\leq r\leq\infty $, $D^{k,r}$ is the homogeneous Sobolev space, which is defined by
\begin{align*}
&D^{k,r}=\big\{u\in L^1_{loc}(\mathbb{R}^3))~\big| ~\|\nabla^ku\|_r< \infty\big\}, ~~~D^k=D^{k,2},\\
&D_0^1=\big\{u\in L^6(\mathbb{R}^3))~\big| ~\|\nabla u\|_2< \infty\big\}.
\end{align*}
For simplicity, let
$$\int f dx=\int_{\mathbb{R}^3}fdx.$$
$\dot{f}:=f_t+u\cdot\nabla f$ denotes the material derivative of $f$.

\begin{definition}
\label{definition 1.1}
Let $T>0$. $(\rho, u, \theta, d)$ is called a strong solution to the compressible
nematic liquid crystal flow (\ref{1.1})-(\ref{1.4}) in $\mathbb{R}^3\times (0, T)$ with initial condition (\ref{1.5}), if for some
${q} \in (3,6]$,
\begin{align*}
&\rho\in C([0,T]; H^1\cap W^{1,{q}}), ~~(u, \theta)\in C([0,T]; D_0^1\cap D^2)\cap L^2(0,T; D^{2, {q}}),\\
&\nabla d\in C([0,T]; H^2)\cap L^2(0,T; H^3), ~~\rho_t \in C([0,T]; L^2\cap L^{{q}}),
~~(u_t, \theta_t) \in L^2(0,T; D_0^1), \\
&(\sqrt{\rho}u_t, \sqrt{\rho}\theta_t)\in L^\infty(0,T; L^2), ~~d_t\in C([0,T]; H^1)\cap L^2(0,T; H^2), ~~|d|=1,
\end{align*}
and $(\rho, u, \theta, d)$ satisfies (\ref{1.1})-(\ref{1.4}) a.e. in $\mathbb{R}^3\times (0, T)$ and fulfills the initial condition (\ref{1.5}).
\end{definition}
\begin{theorem}
\label{theorem 1.1}
Assume that the initial data $\rho_0, u_0, \theta_0$ and $d_0$ satisfy
\begin{align*}
&\rho_0, \theta_0\geq 0, \quad\rho_0\in H^1 \cap W^{1,q},
\quad \sqrt{\rho_0}\theta_0\in L^2,\\
&(u_0, \theta_0)\in D_0^1\cap D^2,\quad \nabla d_0 \in H^2,\quad \text{and } |d_0|=1,
\end{align*}
for some $q \in (3,6]$. Set $\overline{\rho}:=\|\rho_0\|_\infty+1$.
In addition, the
following compatibility condition
\begin{align}
-\mu \Delta  u_0- (\mu+\lambda)\nabla {\rm div}u_0 + \nabla P_0 + \Delta d_0\cdot \nabla d_0=\sqrt{\rho_0}g_1, \label{1.7}\\
\kappa \Delta\theta_0 + \mathcal{Q}(\nabla u_0)+ |\Delta d_0+|\nabla d_0|^2d_0|^2=\sqrt{\rho_0}g_2, \label{1.8}
\end{align}
hold with $g_1$, $g_2 \in L^2$, where $P_0=R\rho_0\theta_0$.

Then, there is a positive constant $\varepsilon_0$ depending only on $R, \gamma, \mu, \lambda$, and $\kappa$,
such that if
$$
\mathcal{N}_0:=\overline{\rho}\big[\|\rho_0\|_3+\overline{\rho}^2(\|\sqrt{\rho_0}u_0\|_2^2+ \|\nabla d_0\|_2^2)\big]
\big[\|\nabla u_0\|_2^2+\overline{\rho}(\|\sqrt{\rho_0}E_0\|_2^2 + \|\nabla^2 d_0\|_2^2)\big]
\leq \varepsilon_0,
$$
where $E_0=\frac{|u_0|^2}{2}+ c_v\theta_0$,
the problem (\ref{1.1})-(\ref{1.5}) has a unique global strong solution.
\end{theorem}

\begin{remark}
\label{remark 1.1}
It is obvious that if there is a initial data $(\rho_0, u_0, \theta_0, d_0)$ satisfying $\mathcal{N}_0 \leq \varepsilon_0$,
then, any $(\rho_{0\lambda}, u_{0\lambda}, \theta_{0\lambda}, d_{0\lambda})$ defined by the scaling transform (\ref{1.6}) with $\lambda\neq 0$
also satisfies $\mathcal{N}_0 \leq \varepsilon_0$.
\end{remark}

\begin{remark}
\label{remark 1.2}
The global well-posedness of strong solutions to the Cauchy
problem for compressible non-isothermal nematic liquid crystal flows with vacuum
as far field density has recently been proved in
\cite{LZ2021}, which needs the small initial data satisfying
 $$\|\rho_0\|_{L^\infty}+ \|\nabla d_0\|_{L^3} \leq \epsilon_0=\epsilon_0(\|\rho_0\|_{1}, \|\sqrt{\rho_0}u_0\|_{2}, \|\nabla u_0\|_{2},
\| \sqrt{\rho_0}E_0\|_{2}, \|\nabla d_0\|_{2}, \|\nabla^2 d_0\|_{2}, \mu, \lambda, R, \gamma, \kappa).$$
Note that the explicit dependence of $\epsilon_0$ on the initial norms was not derived in \cite{LZ2021}.
Therefore, the scaling invariant quantities may not be expected there.
\end{remark}

%The system (\ref{1.1})-(\ref{1.5}) is a combination of the full compressible Navier-Stokes equations
%and the equations for heat flow of harmonic maps with the initial data allowing a vacuum far field or infinite mass. From the study on the global well-posedness of the full compressible Navier-Stokes equations, the %important dissipation estimates $\int^T_0 \|\nabla u\|_2^2dt$ can not be obtained by the basic energy estimates directly as that for the isentropic flows.

Note that the system (\ref{1.1})-(\ref{1.5}) contains the full compressible Navier-Stokes equations as subsystem, it inherits the difficulties of the full compressible Navier-Stokes equations.
A typical difficulty is that the basic energy estimate does not yield the desired dissipation estimates $\int^T_0 \|\nabla u\|_2^2dt$.
Another difficulty is that the presence of the liquid crystal director field $d$
brings strong coupling term $u\cdot \nabla d$ and nonlinear terms $\nabla d\cdot \Delta d$ and $|\nabla d|^2d$.
In order to overcome these difficulties, we adopt the idea in \cite{L2020} to get the $L^\infty(0,T; L^3)$ estimate of $\rho$ and introduce
the spatial $L^2$-norm of $\nabla d$ and $\nabla^2 d$. This motivates us to put smallness assumptions on $\|\rho_0\|^2_\infty\|\sqrt{\rho_0}u_0\|_2\|\sqrt{\rho_0}|u_0|^2\|_2$
and $\|\nabla d_0\|_2\|\nabla^2 d_0\|_2$, which are both scaling invariant. As a result, by continuity arguments, some necessary
lower order time-independent estimates are obtained. Then, we give higher order estimates and eliminate the impact of vacuum by introducing the
effective viscous flux and the material derivative.

The paper is organized as follows. In section 2 and section 3, we derive some lower order and higher order a priori estimates
for the solutions to the Cauchy problem (\ref{1.1})-(\ref{1.5}), respectively.
Section 4 is devoted to proving the global well-posedness.

\section{Time independent lower order a priori estimates}

First, the following local well-posedness of strong solutions with initial vacuum can be established in a similar
way as in \cite{CK2006} and \cite{FLN2018}.
\begin{lemma}
\label{lemma 2.1}(local well-posedness)
Assume that the initial data $(\rho_0\geq 0, u_0, \theta_0\geq 0, d_0)$ satisfies the conditions in Theorem \ref{theorem 1.1}. Then, there is a positive time $T_0$ depending only on $R, c_v, \mu, \lambda, \kappa,$ and $\Phi_0$, such that the system (\ref{1.1})-(\ref{1.5}) admits a unique strong solution in $\mathbb R^3\times(0,T_0)$, where $\Phi_0$ is a positive constant such that
$$
\|\rho_0\|_{W^{1,q}\cap H^1}+\|(u_0,\theta_0)\|_{D_0^1\cap D^2}+\|\nabla d_0\|_{H^2}+\|(\sqrt{\rho_0}\theta_0, g_1, g_2)\|_2\leq \Phi_0.
$$
\end{lemma}

%With this local strong solution $(\rho, u, \theta, d)$ on $\mathbb{R}^3\times [0, T)$ for some $T>0$ at hand, we are ready to show some estimates for this solution.
In the rest of this section, as well as in the next section, we always assume that $(\rho, u, \theta, d)$ is a strong solution to system (\ref{1.1})-(\ref{1.5}) in
$\mathbb{R}^3\times (0, T)$ for some positive time $T$.

\begin{lemma}
\label{lemma 2.2} Assume that $2\mu>\lambda$ . It holds that
\begin{align*}
&\sup_{0\leq t\leq T}(\|\sqrt{\rho}u\|_{2}^2+ \|\nabla d\|_{2}^2)+ \int^T_0(\|\nabla u\|_{2}^2+ \|d_t\|_{2}^2+
\|\nabla^2 d\|_{2}^2)dt\\
&\leq C(\|\sqrt{\rho_0}u_0\|_{2}^2+ \|\nabla d_0\|_{2}^2)+ C\int^T_0 \|\rho\|_3^2\|\nabla\theta\|^2_2dt
+C\int^T_0 \|\nabla d\|_3^2(\|\nabla u\|^2_2+\|\nabla^2 d\|^2_2)dt,
\end{align*}
for a positive constant $C$ depending only on $\mu$ and $\lambda$.
\end{lemma}

{\it\bfseries Proof.}
Multiplying (\ref{1.2}) by $u$, integrating it over $\mathbb{R}^3$, and using integration by parts, one gets
\begin{align}
\label{2.1}
\begin{aligned}
&\frac{1}{2}\frac{d}{dt}\|\sqrt{\rho}u\|_2^2+ \mu\|\nabla u\|_2^2+ (\mu+\lambda)\|{\rm div}u\|_2^2\\
=&-\int \nabla (R\rho\theta)\cdot udx- \int \Delta d\cdot \nabla d\cdot udx\\
\leq& R\|\rho\|_3\|\theta\|_6\|{\rm div}u\|_2+ \|u\|_6\|\nabla d\|_3\|\nabla^2 d\|_2\\
\leq& C\|\rho\|_3\|\nabla\theta\|_2\|{\rm div}u\|_2+ C\|\nabla u\|_2\|\nabla d\|_3\|\nabla^2 d\|_2\\
\leq& (\mu+\lambda)\|{\rm div}u\|_2^2+ C\|\rho\|_3^2\|\nabla\theta\|_2^2+ \frac{\mu}{2}\|\nabla u\|_2^2+ C\|\nabla d\|_3^2\|\nabla^2 d\|_2^2.
\end{aligned}
\end{align}
Using (\ref{1.4}) and the Sobolev inequality, it follows
\begin{align}
\label{2.2}
\begin{aligned}
\frac{d}{dt} \|\nabla d\|_2^2+ \int (|d_t|^2+ |\nabla^2 d|^2)dx=&\int|d_t-\Delta d|^2 dx=\int |u\cdot\nabla d-|\nabla d|^2d|^2dx\\
\leq & C(\|u\|_6^2\|\nabla d\|_3^2+ \|\nabla d\|_4^4) \leq C\|\nabla d\|_3^2(\|\nabla u\|_2^2+ \|\nabla^2 d\|_2^2),
\end{aligned}
\end{align}
where $|\nabla d|^2=-\Delta d\cdot d$ guaranteed by $|d|=1$ was used.
Adding (\ref{2.1}) and (\ref{2.2}) yields
\begin{align*}
&\frac{d}{dt}(\|\sqrt{\rho}u\|_2^2+\|\nabla d\|_2^2)+ \mu\|\nabla u\|_2^2+ \|d_t\|_2^2+ \|\nabla^2 d\|_2^2\\
\leq& C\|\rho\|_3^2\|\nabla\theta\|_2^2+
C\|\nabla d\|_3^2(\|\nabla u\|_2^2+ \|\nabla^2 d\|_2^2).
\end{align*}
The conclusion follows by integrating in $t$ to the above inequality.
\hfill$\Box$

\begin{lemma}
\label{lemma 2.3} It holds that
\begin{align*}
&\sup_{0\leq t\leq T}(\| \nabla^2 d\|_{2}^2+ \|\nabla d\|_{4}^4)+ \int^T_0(\|\nabla d_t\|_{2}^2+ \|\nabla^3 d\|_{2}^2+
\||\nabla d||\nabla^2 d|\|_{2}^2)dt\\
\leq& C(\| \nabla^2 d_0\|_{2}^2+ \|\nabla d_0\|_{4}^4)+C \sup_{0\leq t\leq T}(\|\nabla u\|_2^2\|\nabla d\|_2^2)\int^T_0\|\nabla u\|^6_2dt\\
&+C\sup_{0\leq t\leq T}(\|\nabla d\|_2^3\|\nabla^2 d\|_2^3)\int^T_0\|\nabla^3 d\|^2_2dt,
\end{align*}
for an absolute positive constant $C$.
\end{lemma}

{\it\bfseries Proof.}
Applying the operator $\nabla$ to (\ref{1.4}) yields
\begin{align}
\label{2.3}
\nabla d_t -\nabla\Delta d=-\nabla(u\cdot\nabla d)+|\nabla(\nabla d|^2d),
\end{align}
from which one derives
\begin{align}
\label{2.4}
\begin{aligned}
&\frac{d}{dt} \|\nabla^2 d\|_2^2+ \int (|\nabla d_t|^2+ |\nabla^3 d|^2)dx=
\int|\nabla d_t-\nabla\Delta d|^2 dx=\int |\nabla(u\cdot\nabla d)-|\nabla(\nabla d|^2d)|^2dx\\
\leq & C \int(|\nabla u|^2|\nabla d|^2+ |u|^2|\nabla^2 d|^2+ |\nabla d|^6) dx
+ 2\int |\nabla d|^2|\nabla^2 d|^2dx.
\end{aligned}
\end{align}
Multiplying (\ref{2.3}) by $4|\nabla d|^2\nabla d$, and then integrating it over $\mathbb{R}^3$, one obtains
\begin{align}
\label{2.5}
\begin{aligned}
&\frac{1}{2}\frac{d}{dt}\int|\nabla d|^4dx+ 4\int (|\nabla d|^2|\nabla^2 d|^2+
2|\nabla d|^2|\nabla (|\nabla d|)|^2)dx\\
=& 4\int |\nabla d|^2\nabla d\big(- \nabla(u\cdot\nabla d)+\nabla (|\nabla d|^2d)\big)dx\\
\leq &C \int \big(
|\nabla d|^4|\nabla u|+ |\nabla d|^3|\nabla^2 d||u|+ |\nabla d|^4|\nabla^2 d|
+|\nabla d|^6
\big)dx\\
\leq & \int |\nabla d|^2|\nabla^2 d|^2dx + C \int (|\nabla u|^2|\nabla d|^2+
|u|^2|\nabla^2 d|^2+ |\nabla d|^6) dx,
\end{aligned}
\end{align}
where $|\nabla d|^2=-\Delta d\cdot d$ guaranteed by $|d|=1$ was used.
Adding (\ref{2.4}) and (\ref{2.5}) yields
\begin{align}
\label{2.6}
\begin{aligned}
&\frac{d}{dt} \int(|\nabla^2 d|^2+|\nabla d|^4)dx+ \int (|\nabla d_t|^2+ |\nabla^3 d|^2+|\nabla d|^2|\nabla^2 d|^2)\\
\leq& C \int (|\nabla u|^2|\nabla d|^2+
|u|^2|\nabla^2 d|^2+ |\nabla d|^6) dx.
\end{aligned}
\end{align}
It follows from the Gagliardo-Nirenberg and Young inequalities that
\begin{align}
\label{2.7}
\begin{aligned}
\int |\nabla u|^2|\nabla d|^2 dx \leq& \|\nabla d\|_\infty^2\|\nabla u\|_2^2 \leq C \|\nabla d\|_2^{\frac{1}{2}}\|\nabla^3 d\|_2^{\frac{3}{2}}\|\nabla u\|_2^2\\
\leq& \frac{1}{4}\|\nabla^3 d\|_2^2+ C(\|\nabla d\|_2^2\|\nabla u\|_2^2)\|\nabla u\|_2^6,
\end{aligned}
\end{align}
\begin{align}
\label{2.8}
\begin{aligned}
\int |\nabla d|^6 dx \leq& C\|\nabla d\|_6^2\||\nabla d|^2\|_6\||\nabla d|^2\|_2\leq C\|\nabla d\|_6^2\||\nabla(\nabla d|^2)\|_2\|\nabla d\|_4^2\\
\leq& \frac{1}{2} \||\nabla d||\nabla^2 d|\|_2^2+ C\|\nabla^2 d\|_2^4\|\nabla d\|_4^4
\\
\leq& \frac{1}{2} \||\nabla d||\nabla^2 d|\|_2^2+ C\|\nabla^2 d\|_2^4\|\nabla d\|_3^2\|\nabla^2 d\|_2^2
\\
\leq& \frac{1}{2} \||\nabla d||\nabla^2 d|\|_2^2+ C\|\nabla d\|_2^2\|\nabla^3 d\|_2^2\|\nabla d\|_3^2\|\nabla^2 d\|_2^2
\\
\leq& \frac{1}{2} \||\nabla d||\nabla^2 d|\|_2^2+ C\|\nabla d\|_2^3\|\nabla^2 d\|_2^3\|\nabla^3 d\|_2^2,
\end{aligned}
\end{align}
and
\begin{align}
\label{2.9}
\begin{aligned}
\int |u|^2|\nabla^2 d|^2 dx \leq& C\|u\|_6^2\|\nabla^2 d\|_6\|\nabla^2 d\|_2\leq
C\|\nabla u\|_2^2\|\nabla d\|_2^{\frac{1}{2}}\|\nabla^3 d\|_2^{\frac{3}{2}}\\
\leq& \frac{1}{4}\|\nabla^3 d\|_2^2+ C(\|\nabla d\|_2^2\|\nabla u\|_2^2)\|\nabla u\|_2^6.
\end{aligned}
\end{align}
Substituting (\ref{2.7})-(\ref{2.9}) into (\ref{2.6}) yields
\begin{align*}
&\frac{d}{dt} \int(|\nabla^2 d|^2+|\nabla d|^4)dx+ \frac{1}{2}\int (|\nabla d_t|^2+ |\nabla^3 d|^2+|\nabla d|^2|\nabla^2 d|^2)\\
\leq&  C\|\nabla d\|_2^3\|\nabla^2 d\|_2^3\|\nabla^3 d\|_2^2 +C(\|\nabla d\|_2^2\|\nabla u\|_2^2)\|\nabla u\|_2^6,
\end{align*}
which implies the conclusion by integrating in $t$.
\hfill$\Box$

\begin{lemma}
\label{lemma 2.4} It holds that
\begin{align*}
&\sup_{0\leq t\leq T}\|\sqrt{\rho}E\|_{2}^2+ \int^T_0(\|\nabla \theta\|_{2}^2+ \||u||\nabla u|\|_{2}^2)dt\\
\leq& C\|\sqrt{\rho_0}E_0\|_{2}^2+C \int^T_0\|\rho\|_\infty\|\rho\|_3^\frac{1}{2}\|\sqrt{\rho}\theta\|_2(\|\nabla \theta\|^2_2+ \||u||\nabla u|\|_{2}^2)dt+ C\sup_{0\leq t\leq T}\|\nabla^2 d\|_2^4\int^T_0\|\nabla u\|^2_2dt\\
&
+ C\sup_{0\leq t\leq T}(\|\nabla d\|_2^2\|\nabla^2 d\|_2^2)
\int^T_0\||\nabla d||\nabla^2 d|\|_{2}^2dt +C\sup_{0\leq t\leq T}(\|\nabla d\|_2\|\nabla^2 d\|_2)\int^T_0\|\nabla^3 d\|^2_2dt,
\end{align*}
for a positive constant $C$ depending only on $R, c_v, \mu, \lambda$, and $\kappa$, where $E=\frac{|u|^2}{2}+c_v\theta$.
\end{lemma}

{\it\bfseries Proof.}
From the energy $E$, multiplying (\ref{1.2}) by $u$ and adding the resultant to (\ref{1.3}), one has
\begin{align}
\label{2.10}
\rho(E_t+u\cdot\nabla E)+{\rm div}(uP)-\kappa\Delta\theta={\rm div}(\mathcal{S}\cdot u)-\Delta d\cdot\nabla d\cdot u+|\Delta d+|\nabla d|^2d|^2,
\end{align}
where $\mathcal{S}=\mu(\nabla u+(\nabla u)^t)+ \lambda{\rm div} uI$. Then, multiplying (\ref{2.10}) with
$E$ and integrating it over $\mathbb{R}^3$ yield
\begin{align*}
\frac{1}{2}\frac{d}{dt}\|\sqrt{\rho}E\|_2^2+\kappa c_v\|\nabla\theta\|_2^2 \leq & \frac{\kappa c_v}{2}\|\nabla\theta\|_2^2 +
C\||u||\nabla u|\|_2^2 + C\int \rho^2\theta^2|u|^2dx\\
& - \int (\Delta d\cdot\nabla d\cdot u) Edx +\int |\Delta d+|\nabla d|^2d|^2Edx
\end{align*}
and, thus,
\begin{align}
\label{2.11}
\begin{aligned}
\frac{d}{dt}\|\sqrt{\rho}E\|_2^2+\kappa c_v\|\nabla\theta\|_2^2 \leq &
C\||u||\nabla u|\|_2^2 + C\int \rho^2\theta^2|u|^2dx\\
& - 2\int (\Delta d\cdot\nabla d\cdot u) Edx +2\int |\Delta d+|\nabla d|^2d|^2Edx.
\end{aligned}
\end{align}
One can rewrite the right hand side of (\ref{1.2}) in divergence form since
\begin{align*}
-\Delta d\cdot\nabla d=-{\rm div}\bigg(\nabla d \odot \nabla d -\frac{1}{2}|\nabla d|^2\mathbb{I}_3\bigg),
\end{align*}
where
$
\nabla d \odot \nabla d \triangleq (d_{x_i}\cdot d_{x_j})_{3 \times 3}
$
and $\mathbb{I}_3$ denotes the identity matrix of order 3.
By the Sobolev inequality and integration by parts, one deduces
\begin{align*}
&-\int (\Delta d\cdot\nabla d\cdot u) Edx\\
=&-\int {\rm div}\bigg(\nabla d \odot \nabla d -\frac{1}{2}|\nabla d|^2\mathbb{I}_3\bigg)\cdot u Edx\\
=&\int \bigg(\nabla d \odot \nabla d -\frac{1}{2}|\nabla d|^2\mathbb{I}_3\bigg):\nabla(uE)dx\\
\leq& C\|\nabla d\|_6^2\|\nabla u\|_2\|E\|_6+ \frac{\kappa c_v}{16}\|\nabla \theta\|_2^2 +C\||u|\nabla u\|_2^2+C\int|\nabla d|^4|u|^2dx\\
\leq& \frac{\kappa c_v}{8}\|\nabla \theta\|_2^2 +C\||u|\nabla u\|_2^2+ C\|\nabla^2 d\|_2^4\|\nabla u\|_2^2
+C\||u|^2\|_6\|\nabla d\|_2\|\nabla d\|_6\||\nabla d|^2\|_6\\
\leq& \frac{\kappa c_v}{8}\|\nabla \theta\|_2^2 +C\||u|\nabla u\|_2^2+ C\|\nabla^2 d\|_2^4\|\nabla u\|_2^2
+C\|\nabla d\|_2^2\|\nabla^2 d\|_2^2\||\nabla d||\nabla^2 d|\|_2^2,
\end{align*}
and
\begin{align*}
&\int |\Delta d+|\nabla d|^2d|^2E
\leq  2\int |\Delta d|^2E dx+2\int |\nabla d|^4E dx\\
\leq & 2\int |\nabla d||\nabla^3 d|E dx+ 2\int |\nabla d||\nabla^2 d||\nabla E| dx+ \int|\nabla d|^4|u|^2dx+2c_v\int|\nabla d|^4\theta dx\\
\leq & C\|\nabla d\|_3\|\nabla^3 d\|_2\|E\|_6+ C\|\nabla d\|_3\|\nabla^2 d\|_6\|\nabla E\|_2+
 C(\||u|^2\|_6+c_v\|\theta\|_6)\|\nabla d\|_2\|\nabla d\|_6\||\nabla d|^2\|_6\\
\leq & C (\||u||\nabla u|\|_2+c_v\|\nabla\theta\|_2)(\|\nabla d\|_3\|\nabla^3 d\|_2+\|\nabla d\|_2\|\nabla^2 d\|_2\||\nabla d||\nabla^2 d|\|_2)\\
\leq& \frac{\kappa c_v}{8}\|\nabla \theta\|_2^2 +C\||u|\nabla u\|_2^2+ C\|\nabla d\|_3^2\|\nabla^3 d\|_2^2
+C\|\nabla d\|_2^2\|\nabla^2 d\|_2^2\||\nabla d||\nabla^2 d|\|_2^2.
\end{align*}
Putting the above two inequalities into (\ref{2.11}) leads to
\begin{align}
\label{2.12}
\begin{aligned}
\frac{d}{dt}\|\sqrt{\rho}E\|_2^2+\frac{\kappa c_v}{2}\|\nabla\theta\|_2^2 \leq &
C\||u||\nabla u|\|_2^2 + C\int \rho^2\theta^2|u|^2dx\\
&+ C\|\nabla^2 d\|_2^4\|\nabla u\|_2^2+ C\|\nabla d\|_3^2\|\nabla^3 d\|_2^2
+C\|\nabla d\|_2^2\|\nabla^2 d\|_2^2\||\nabla d||\nabla^2 d|\|_2^2.
\end{aligned}
\end{align}

To control the term $\||u||\nabla u|\|_2^2$ in (\ref{2.12}), we need to multiply (\ref{1.2}) with
$|u|^2u$ to obtain that
\begin{align}
\label{2.13}
\begin{aligned}
\frac{1}{4}\frac{d}{dt}&\|\sqrt{\rho}|u|^2\|_2^2-\int(\mu\Delta u+ (\mu+\lambda)\nabla{\rm div}u)\cdot|u|^2udx\\
&=-\int P {\rm div}(|u|^2u)dx- \int \bigg(\nabla d \odot \nabla d -\frac{1}{2}|\nabla d|^2\mathbb{I}_3\bigg) {\rm div}(|u|^2u)dx\\
& \leq C\int  (\rho\theta+ |\nabla d|^2)|\nabla u||u|^2dx\\
& \leq \frac{1}{2}\left(\mu-\frac{\lambda}{2}\right)\||u||\nabla u|\|_2^2+ C\int  (\rho^2\theta^2|u|^2+ |\nabla d|^4|u|^2)dx\\
& \leq \frac{1}{2}\left(\mu-\frac{\lambda}{2}\right)\||u||\nabla u|\|_2^2+ C\int  \rho^2\theta^2|u|^2dx
+C\||u|^2\|_6\|\nabla d\|_2\|\nabla d\|_6\||\nabla d|^2\|_6\\
&\leq \left(\mu-\frac{\lambda}{2}\right)\||u||\nabla u|\|_2^2+ C\int \rho^2\theta^2|u|^2dx
+C\|\nabla d\|_2^2\|\nabla^2 d\|_2^2\||\nabla d||\nabla^2 d|\|_2^2.
\end{aligned}
\end{align}
By direct computation, one has
\begin{align*}
-\int(\mu\Delta u+ (\mu+\lambda)\nabla{\rm div}u)\cdot|u|^2udx \geq (2\mu-\lambda)\||u||\nabla u|\|_2^2.
\end{align*}
Hence, it follows from the above inequality and (\ref{2.13}) that
\begin{align}
\label{2.14}
\frac{d}{dt}\|\sqrt{\rho}|u|^2\|_2^2 +2(2\mu-\lambda)\||u||\nabla u|\|_2^2
\leq  C\int \rho^2\theta^2|u|^2dx+C\|\nabla d\|_2^2\|\nabla^2 d\|_2^2\||\nabla d||\nabla^2 d|\|_2^2.
\end{align}

Now, multiplying (\ref{2.14}) by $N>0$, which is a sufficiently large number and depending only on $R, c_v, \mu, \lambda$, and $\kappa$, then adding the resultant to (\ref{2.12}), one gets
\begin{align}
\label{2.15}
\begin{aligned}
&\frac{d}{dt}(\|\sqrt{\rho}E\|_2^2 + N\|\sqrt{\rho}|u|^2\|_2^2)  +\frac{\kappa c_v}{2}\|\nabla\theta\|_2^2 + (2\mu-\lambda)N\||u||\nabla u|\|_2^2\\
\leq &C\|\rho\|_\infty\|\rho\|_3^{\frac{1}{2}}\|\sqrt{\rho}\theta\|_2\|\nabla\theta\|_2\||u||\nabla u|\|_2
+ C\|\nabla^2 d\|_2^4\|\nabla u\|_2^2+ C\|\nabla d\|_2\|\nabla^2 d\|_2 \|\nabla^3 d\|_2^2\\
&+C\|\nabla d\|_2^2\|\nabla^2 d\|_2^2\||\nabla d||\nabla^2 d|\|_2^2,
\end{aligned}
\end{align}
where we have used the fact
\begin{align}
\label{2.16}
\int \rho^2\theta^2|u|^2dx\leq C\|\sqrt{\rho}\theta\|_2\|\theta\|_6\||u|^2\|_6\|\rho\|_9^{\frac{3}{2}}
\leq
C\|\rho\|_\infty\|\rho\|_3^{\frac{1}{2}}\|\sqrt{\rho}\theta\|_2\|\nabla\theta\|_2\||u||\nabla u|\|_2.
\end{align}
The proof can be completed by integrating (\ref{2.15}) over $[0, t]$.
\hfill$\Box$

\begin{lemma}
\label{lemma 2.5} It holds that
\begin{align*}
\sup_{0\leq t\leq T}\|\rho\|_{3}^3+ \int^T_0\int \rho^3Pdxdt
\leq& C\|\rho_0\|_3^3+ C\sup_{0\leq t\leq T}(\|\rho\|_\infty^{\frac{2}{3}}\|\sqrt{\rho}u\|_2^{\frac{1}{3}}\|\sqrt{\rho}|u|^2\|_2^{\frac{1}{3}}
\|\rho\|_3^3)\\
& +C\int^T_0(\|\rho\|_\infty^2\|\rho\|_3^2\|\nabla u\|^2_2)dt
%+ C\int^T_0\int \rho^3|\nabla d|^2dxdt.
+C\int^T_0(\|\rho\|_\infty\|\rho\|_3^2\|\nabla^2 d\|^2_2)dt,
\end{align*}
for a positive constant $C$ depending only on $R, c_v, \mu, \lambda$, and $\kappa$.
\end{lemma}

{\it\bfseries Proof.}
If we apply the operator $\Delta^{-1}{\rm div}$
to (\ref{1.2}), it holds that
\begin{align}
\label{2.17}
\Delta^{-1}{\rm div}(\rho u)_t+ \Delta^{-1}{\rm div}{\rm div}(\rho u \otimes u)-(2\mu+\lambda){\rm div}u+P
=-\Delta^{-1}{\rm div}{\rm div}\big(\nabla d \odot \nabla d -\frac{1}{2}|\nabla d|^2\mathbb{I}_3\big).
\end{align}
Taking advantage of (\ref{1.1}), we have
\begin{align}
\label{2.17-1}
(\rho^3)_t+ {\rm div}(u\rho^3)+2{\rm div}u\rho^3=0.
\end{align}
Then, multiplying (\ref{2.17}) by $\rho^3$ and using above equality, we obtain
\begin{align}
\label{2.18}
\begin{aligned}
\frac{2\mu+\lambda}{2}\big( (\rho^3)_t+{\rm div}(u\rho^3)\big)+& \rho^3 P+\rho^3 \Delta^{-1}{\rm div}(\rho u)_t+ \rho^3\Delta^{-1}{\rm div}{\rm div}(\rho u \otimes u)\\
=&-\rho^3\Delta^{-1}{\rm div}{\rm div}\big(\nabla d \odot \nabla d -\frac{1}{2}|\nabla d|^2\mathbb{I}_3\big).
\end{aligned}
\end{align}
Using (\ref{2.17-1}), it follows
\begin{align*}
\int \rho^3 \Delta^{-1}{\rm div}(\rho u)_tdx=& \frac{d}{dt}\int \rho^3 \Delta^{-1}{\rm div}(\rho u) dx + \int [{\rm div}(\rho^3 u)+ 2{\rm div}u\rho^3]\Delta^{-1}{\rm div}(\rho u)dx\\
=& \int [2{\rm div}u\rho^3\Delta^{-1}{\rm div}(\rho u)- \rho^3 u\cdot\nabla\Delta^{-1}{\rm div}(\rho u)]dx + \frac{d}{dt}\int \rho^3 \Delta^{-1}{\rm div}(\rho u) dx.
\end{align*}
Thanks to this, it follows from integrating (\ref{2.18}) over $\mathbb{R}^3$ that
\begin{align*}
\frac{d}{dt}\int& \left(\frac{2\mu+\lambda}{2}+ \Delta^{-1}{\rm div}(\rho u)\right)\rho^3dx+\int\rho^3Pdx\\
=&
 \int\left[\rho^3\big(u\cdot \nabla\Delta^{-1}{\rm div}(\rho u)-\Delta^{-1}{\rm div}{\rm div}(\rho u \otimes u)\big)-2{\rm div}u\rho^3\Delta^{-1}{\rm div}(\rho u)\right]dx\\
&- \int \rho^3\Delta^{-1}{\rm div}{\rm div}\bigg(\nabla d \odot \nabla d -\frac{1}{2}|\nabla d|^2\mathbb{I}_3\bigg)dx.
\end{align*}
The conclusion in this lemma then follows from the same estimates as in Proposition 2.4 in \cite{L2020} and the following bound for the last term in above equality
\begin{align*}
&\int \rho^3\left|\Delta^{-1}{\rm div}{\rm div}\bigg(\nabla d \odot \nabla d -\frac{1}{2}|\nabla d|^2\mathbb{I}_3\bigg)\right|dx\\
 \leq& C\|\rho\|_\infty\|\rho\|_3^2\bigg\|\Delta^{-1}{\rm div}{\rm div}\bigg(\nabla d \odot \nabla d -\frac{1}{2}|\nabla d|^2\mathbb{I}_3\bigg)\bigg\|_3\\
\leq& C\|\rho\|_\infty\|\rho\|_3^2\| \nabla d\|^2_6
\leq C\|\rho\|_\infty\|\rho\|_3^2\|\nabla^2 d\|^2_2,
\end{align*}
where the elliptic estimates were applied.
\hfill$\Box$

\begin{lemma}
\label{lemma 2.6}
Assume that
\begin{align*}
\sup_{0\leq t\leq T}\|\rho\|_\infty \leq 4\overline{\rho}.
\end{align*}
Then, it holds that
\begin{align*}
&\sup_{0\leq t\leq T}\|\nabla u\|_{2}^2+ \int^T_0 \bigg\|\bigg(\sqrt{\rho}u_t, \frac{\nabla G}{\sqrt{\overline{\rho}}}, \frac{\nabla \omega}{\sqrt{\overline{\rho}}}\bigg)\bigg\|_2^2dt\\
\leq& C\|\nabla u_0\|_2^2+ C\overline{\rho}\sup_{0\leq t\leq T}\|\sqrt{\rho}\theta\|_2^2
+C\overline{\rho}^3\int^T_0\|\nabla u\|_2^4(\|\nabla u\|_2^2+\overline{\rho}\|\sqrt{\rho}\theta\|_2^2 )dt\\
&+C\int^T_0 (\overline{\rho} + \overline{\rho}^2\|\rho\|_3^\frac{1}{2}\|\sqrt{\rho}\theta\|_2)(\|\nabla\theta\|_2^2+\||u||\nabla u|\|_2^2)dt\\
&+C\|\nabla d_0\|_4^4+ C\|\nabla d\|_4^4+ \varepsilon_1\overline{\rho}\int^T_0\|\nabla\theta\|_2^2dt +
C\overline{\rho}\sup_{0\leq t\leq T}(\|\nabla d\|_{2}\|\nabla^2 d\|_{2})\int^T_0(\|\nabla d_t\|_2^2+\|\nabla^3 d\|_2^2)dt\\
&+C\overline{\rho}\sup_{0\leq t\leq T}(\|\nabla d\|_{2}^2\|\nabla^2 d\|_{2}^2)\int^T_0\|\nabla (|\nabla d|^2) \|_2^2dt,
\end{align*}
where $G=(2\mu+\lambda){\rm div}u-p$, $\omega=\nabla\times u$ and the constant $C>0$ depending only on $R, c_v, \mu, \lambda$, and $\kappa$.
\end{lemma}

{\it\bfseries Proof.}
Multiplying (\ref{1.2}) by $u_t$, it follows from integration by parts that
\begin{align}
\label{2.19}
\begin{aligned}
 \frac{1}{2}\frac{d}{dt}&\big(\mu\|\nabla u\|_2^2+ (\mu+\lambda)\|{\rm div}u\|_2^2\big)-\int P {\rm div}u_tdx+ \|\sqrt{\rho}u_t\|_2^2\\
 &=-\int \rho(u\cdot\nabla)u\cdot u_tdx-\int \Delta d\cdot\nabla d\cdot u_tdx.
\end{aligned}
\end{align}
By the definition of effective viscous flux $G$, it is easy to see ${\rm div}u=\frac{G+P}{2\mu+\lambda}$, which implies
\begin{align}
\label{2.20}
\begin{aligned}
-\int P {\rm div}u_tdx=&-\frac{d}{dt}\int P{\rm div}udx+ \int P_t{\rm div}udx\\
=& -\frac{d}{dt}\int P{\rm div}udx+\frac{1}{2(2\mu+\lambda)}\frac{d}{dt}\|P\|_2^2+ \frac{1}{2\mu+\lambda}\int P_tGdx.
\end{aligned}
\end{align}
On the other hand, it follows from (\ref{1.3}) that
\begin{align*}
P_t=(\gamma-1)\big(\mathcal{Q}(\nabla u)-P{\rm div}u+ \kappa\Delta\theta + |\Delta d+|\nabla d|^2d|^2\big)-{\rm div}(uP),
\end{align*}
which leads to
\begin{align}
\label{2.21}
\begin{aligned}
\int P_tGdx= \int\left[(\gamma-1)\big(\mathcal{Q}(\nabla u)-P{\rm div}u+ |\Delta d+|\nabla d|^2d|^2\big)G +
\big(uP-\kappa(\gamma-1)\nabla\theta\big)\cdot\nabla G\right]dx.
\end{aligned}
\end{align}
Due to $\|\nabla u\|_2^2=\|\omega\|_2^2+\|{\rm div}u\|_2^2$ and combining with (\ref{2.19})-(\ref{2.21}), one can deduce that
\begin{align}
\label{2.22}
\begin{aligned}
&\frac{1}{2}\frac{d}{dt}\left(\mu\|\omega\|_2^2+ \frac{\|G\|_2^2}{2\mu+\lambda}\right)+ \|\sqrt{\rho}u_t\|_2^2\\
=&
-\int \rho(u\cdot \nabla)u\cdot u_tdx-\int\Delta d \cdot \nabla d\cdot u_tdx
+\frac{1}{2\mu+\lambda}\int \big(\kappa(\gamma-1)\nabla\theta-uP\big)\cdot\nabla G dx\\
&- \frac{\gamma-1}{2\mu+\lambda}
\int \big(\mathcal{Q}(\nabla u)-p{\rm div}u+ |\Delta d+|\nabla d|^2d|^2\big)Gdx.
\end{aligned}
\end{align}

In order to bound the right hand side of (\ref{2.22}),
we need to reformulate (\ref{1.2}) in the following form with the help of $\Delta u= \nabla {\rm div}u- \nabla\times\nabla\times u$:
\begin{align}
\label{2.23}
\rho(u_t+u\cdot \nabla u)=\nabla G-\mu\nabla\times \omega-\Delta d\cdot\nabla d.
\end{align}
Then, multiplying both sides of (\ref{2.23}) by $\nabla G$, it follows
\begin{align*}
\|\nabla G\|_2^2 =& \int \big(\rho(u_t+u\cdot \nabla u) \cdot\nabla G + \Delta d\cdot\nabla d\cdot\nabla G \big)dx\\
 \leq &\int \left( \frac{|\nabla G|^2}{2}+2\overline{\rho}\rho|u_t|^2 \right)dx +\int \big(\rho (u\cdot \nabla) u\cdot\nabla G+ \Delta d\cdot\nabla d\cdot\nabla G \big)dx
\end{align*}
where $\int \nabla G \cdot \nabla \times\omega dx=0$ and $\|\rho\|_\infty\leq 4\overline{\rho}$ were used.
This gives that
\begin{align}
\label{2.24}
\frac{\|\nabla G\|_2^2}{16\overline{\rho}}\leq \frac{1}{4}\|\sqrt{\rho}u_t\|_2^2
+ \frac{1}{8\overline{\rho}}\int \big(\rho (u\cdot \nabla) u\cdot\nabla G+ \Delta d\cdot\nabla d\cdot\nabla G \big)dx.
\end{align}
Similarly, one has
\begin{align}
\label{2.25}
\frac{\mu^2\|\nabla \omega\|_2^2}{16\overline{\rho}}\leq \frac{1}{4}\|\sqrt{\rho}u_t\|_2^2
+ \frac{1}{8\overline{\rho}}\int \big(\rho (u\cdot \nabla) u\cdot\nabla \omega+ \Delta d\cdot\nabla d\cdot\nabla \omega \big)dx.
\end{align}
Putting (\ref{2.24}) and (\ref{2.25}) into (\ref{2.22}), one gets
\begin{align}
\label{2.26}
\begin{aligned}
&\frac{1}{2}\frac{d}{dt}\left(\mu\|\omega\|_2^2+ \frac{\|G\|_2^2}{2\mu+\lambda}\right)+ \frac{1}{2}\|\sqrt{\rho}u_t\|_2^2
+\frac{1}{16\overline{\rho}}(\|\nabla G\|_2^2+ \mu^2\|\nabla \omega\|_2^2)\\
\leq& C\int \rho|u||\nabla u|\left(|u_t|^2+\frac{1}{\overline{\rho}}(|\nabla G|+|\nabla\omega|) \right)dx+
C\int (|\nabla\theta|+\rho\theta|u|)|\nabla G|dx\\
& + C\int (|\nabla u|^2+\rho\theta|\nabla u|)|G|dx
-\int\Delta d\cdot\nabla d\cdot u_tdx
+C\int \big(|\Delta d+|\nabla d|^2d|^2\big)Gdx\\
&+ \frac{C}{\overline{\rho}}\int |\Delta d||\nabla d|(|\nabla G|+ |\nabla\omega|) dx=:\sum^6_{i=1}I_i.
\end{aligned}
\end{align}
Estimate $I_i, i = 1, 2, \ldots, 6$ are given as follows. It follows from the H\"{o}lder and Young inequalities that
\begin{align*}
I_1\leq& C\sqrt{\overline{\rho}}\||u||\nabla u|\|_2\|\sqrt{\rho}u_t\|_2+C\||u||\nabla u|\|_2(\|\nabla G\|_2+ \|\nabla\omega\|_2)\\
\leq& \frac{1}{12}\|\sqrt{\rho}u_t\|_2^2 + \frac{1}{192\overline{\rho}}(\|\nabla G\|_2^2+ \mu^2\|\nabla\omega\|_2^2)+  C\overline{\rho}\||u||\nabla u|\|_2^2,
\end{align*}
and
\begin{align*}
I_2\leq& C\|\nabla \theta\|_2\|\nabla G\|_2+ \|\rho\theta u\|_2\|\nabla G\|_2\\
\leq& C\|\nabla \theta\|_2\|\nabla G\|_2+ C \sqrt{\overline{\rho}}\|\rho\|_3^{\frac{1}{4}}\|\sqrt{\rho}\theta\|_2^{\frac{1}{2}}\|\nabla\theta\|_2^{\frac{1}{2}}\||u||\nabla u|\|_2^{\frac{1}{2}}
 \|\nabla G\|_2\\
\leq& \frac{1}{192\overline{\rho}}\|\nabla G\|_2^2 + C\big(\overline{\rho}^2\|\rho\|_3^{\frac{1}{2}}\|\sqrt{\rho}\theta\|_2+\overline{\rho}
\big)(\|\nabla\theta\|_2^2+ \||u||\nabla u|\|_2^2),
\end{align*}
where (\ref{2.16}) was used in $I_2$. For $I_3$, noticing that
\begin{align}
\label{2.27}
\|\nabla u\|_6\leq C(\|\omega\|_6+\|{\rm div}u\|_6)\leq C(\|\omega\|_6+\|G\|_6+ \|\rho\theta\|_6)
\leq C(\|\nabla\omega\|_2+\|\nabla G\|_2+ \overline{\rho}\|\nabla\theta\|_2),
\end{align}
it follows from the H\"{o}lder, Young and Sobolev inequalities that
\begin{align*}
I_3\leq& C\|\nabla u\|_2\|\nabla u\|_6\| G\|_3+ \|\nabla u\|_2\|\rho\theta\|_6\| G\|_3\\
\leq& C\|\nabla u\|_2(\|\nabla G\|_2+ \|\nabla \omega\|_2 +  \overline{\rho}\|\nabla\theta\|_2)\|G\|_2^{\frac{1}{2}}\|\nabla G\|_2^{\frac{1}{2}}+ C\overline{\rho}\|\nabla u\|_2\|\nabla \theta\|_2\|G\|_2^{\frac{1}{2}}\|\nabla G\|_2^{\frac{1}{2}}\\
\leq& \frac{1}{192\overline{\rho}}(\|\nabla G\|_2^2 +\mu^2\|\nabla \omega\|_2^2)+
 C \overline{\rho}^3 \|\nabla u\|_2^4\| G\|_2^2+ C \overline{\rho}\|\nabla \theta\|_2^2.
\end{align*}
Note that
\begin{align*}
I_4=-\int {\rm div}\bigg(\nabla d \odot \nabla d -\frac{1}{2}|\nabla d|^2\mathbb{I}_3\bigg)\cdot u_t dx
=\int (\nabla d \odot \nabla d) :\nabla u_t dx - \frac{1}{2}\int|\nabla d|^2 {\rm div} u_t dx.
\end{align*}
Using the Young and Sobolev inequalities, one obtains
\begin{align*}
&\int (\nabla d \odot \nabla d) :\nabla u_t dx\\
=& \frac{d}{dt}\int (\nabla d \odot \nabla d) :\nabla u dx
- \int (\nabla d_t \odot \nabla d) :\nabla u dx-\int (\nabla d \odot \nabla d_t) :\nabla u dx\\
\leq & \frac{d}{dt}\int (\nabla d \odot \nabla d) :\nabla u dx
+C \|\nabla d\|_3\|\nabla d_t\|_2\|\nabla u\|_6\\
\leq & \frac{d}{dt}\int (\nabla d \odot \nabla d):\nabla u dx
+\frac{\varepsilon}{4\overline{\rho}}\|\nabla u\|_6^2 +C\overline{\rho} \|\nabla d\|_3^2\|\nabla d_t\|_2^2,
\end{align*}
where $\varepsilon>0$ is a sufficiently small constant.
Similarly,
\begin{align*}
\int |\nabla d|^2{\rm div} u_t dx
\leq & \frac{d}{dt}\int  |\nabla d|^2 {\rm div}u dx
+\frac{\varepsilon}{4\overline{\rho}}\|\nabla u\|_6^2 +C\overline{\rho} \|\nabla d\|_3^2\|\nabla d_t\|_2^2.
\end{align*}
Therefore, we get
\begin{align*}
I_4
\leq  \frac{d}{dt}\int [(\nabla d \odot \nabla d) :\nabla u + |\nabla d|^2 {\rm div}u] dx
+\frac{\varepsilon}{2\overline{\rho}}\|\nabla u\|_6^2 +C\overline{\rho} \|\nabla d\|_3^2\|\nabla d_t\|_2^2.
\end{align*}
Using (\ref{2.27}) again, we obtain
\begin{align*}
I_4
\leq&  \frac{d}{dt}\int [(\nabla d \odot \nabla d) :\nabla u + |\nabla d|^2 {\rm div}u] dx
\\&+
\frac{1}{192\overline{\rho}}(\|\nabla\omega\|_2^2+\|\nabla G\|_2^2)+\frac{\varepsilon_1\overline{\rho}}{2} \|\nabla\theta\|_2^2 +C\overline{\rho} \|\nabla d\|_3^2\|\nabla d_t\|_2^2,
\end{align*}
where $\varepsilon_1=\frac{C\varepsilon}{2}$ is sufficiently small.

Now, let us turn to $I_5$ and $I_6$. By virtue of the H\"{o}lder, Young and Sobolev inequalities, one deduces
\begin{align*}
&I_5 = C\int \big(|\Delta d+|\nabla d|^2d|^2\big)Gdx\\
 \leq& C \int (|\nabla d||\nabla^3 d||G|+ |\nabla d||\nabla^2 d||\nabla G|)dx
+ C \|\nabla d\|_2\|\nabla d\|_6\||\nabla d|^2\|_6\| G\|_6\\
\leq& C \|\nabla d\|_3\|\nabla^3 d\|_2\|\nabla G\|_2
+ C \|\nabla d\|_2\|\nabla^2 d\|_2\|\nabla(|\nabla d|^2)\|_2\|\nabla G\|_2\\
\leq& \frac{1}{192\overline{\rho}}\|\nabla G\|_2^2 + C\overline{\rho}\|\nabla d\|_3^2\|\nabla^3 d\|_2^2+  C\overline{\rho} \|\nabla d\|_2^2\|\nabla^2 d\|_2^2\|\nabla(|\nabla d|^2)\|_2^2,
\end{align*}
where $|\Delta d+ |\nabla d|^2d|^2= |\Delta d|^2+ 2\Delta d\cdot d|\nabla d|^2+ |\nabla d|^4 = |\Delta d|^2-|\nabla d|^4$ was used,
and
\begin{align*}
I_6 \leq \frac{C}{4\overline{\rho}} \|\Delta d\|_6\|\nabla d\|_3 (\|\nabla G\|_2 +\|\nabla \omega\|_2)
\leq \frac{1}{192\overline{\rho}}(\|\nabla G\|_2^2  + \mu^2\|\nabla \omega\|_2^2)+ C\|\nabla d\|_3^2\|\nabla^3 d\|_2^2.
\end{align*}
Putting all these estimates for $I_i, i=1,2,\ldots, 6$ into (\ref{2.26}) leads to
\begin{align}
\label{2.28}
\begin{aligned}
&\frac{d}{dt}\left(\mu\|\omega\|_2^2+ \frac{\|G\|_2^2}{2\mu+\lambda}\right)+ \frac{1}{2}\|\sqrt{\rho}u_t\|_2^2
+\frac{1}{16\overline{\rho}}(\|\nabla G\|_2^2+ \mu^2\|\nabla \omega\|_2^2)\\
\leq& 2\frac{d}{dt}\int [(\nabla d \odot \nabla d) :\nabla u + |\nabla d|^2 {\rm div}u] dx+C (\overline{\rho} + \overline{\rho}^2\|\rho\|_3^\frac{1}{2}\|\sqrt{\rho}\theta\|_2)(\|\nabla\theta\|_2^2+\||u||\nabla u|\|_2^2)\\
&+ C\overline{\rho}^3\|\nabla u\|_2^4\|G\|_2^2
+\varepsilon_1 \overline{\rho}\|\nabla\theta\|_2^2
+C\overline{\rho} \|\nabla d\|_3^2(\|\nabla d_t\|_2^2
+\|\nabla^3 d\|_2^2)\\
&+  C\overline{\rho} \|\nabla d\|_2^2\|\nabla^2 d\|_2^2\|\nabla(|\nabla d|^2)\|_2^2.
\end{aligned}
\end{align}
Note that
\begin{align}
\label{2.29}
\|\nabla u\|_2 \leq C(\|\omega\|_2+\|G\|_2+\|\rho\theta\|_2) \leq C(\|\omega\|_2+\|G\|_2+\sqrt{\overline{\rho}}\|\sqrt{\rho}\theta\|_2)
\end{align}
and
\begin{align}
\label{2.30}
\int [ (\nabla d \odot \nabla d) :\nabla u + |\nabla d|^2 {\rm div}u] dx \leq C\|\nabla u\|_2\||\nabla d|^2\|_2 \leq
\varepsilon_2\|\nabla u\|_2^2+ C\|\nabla d\|_4^4.
\end{align}
Integrating (\ref{2.28}) in $t$, substituting (\ref{2.29}) and (\ref{2.30}) into the resultant and choosing $\varepsilon_2$ small enough, one gets the desired result.
\hfill$\Box$

\begin{lemma}
\label{lemma 2.7}
Assume that
\begin{align*}
\sup_{0\leq t\leq T}\|\rho\|_\infty \leq 4\overline{\rho}.
\end{align*}
Then, it holds that
\begin{align*}
&\sup_{0\leq t\leq T}\|\rho\|_\infty\\
\leq& \|\rho_0\|_\infty e^{C\overline{\rho}^{\frac{2}{3}}\sup_{0\leq t\leq T}\|\sqrt{\rho}u\|_2^{\frac{1}{3}}\|\sqrt{\rho}|u|^2\|_2^{\frac{1}{3}} +
C\overline{\rho}\int^T_0\|\nabla u\|_2\|(\nabla G, \nabla \omega, \overline{\rho}\nabla\theta)\|_2dt
+ C\left(\int^T_0 \|\nabla^2 d\|_2dt\int^T_0\|\nabla^3 d\|_2dt\right)^\frac{1}{2}},
\end{align*}
for a positive constant $C$ depending only on $R, c_v, \mu, \lambda$, and $\kappa$.
\end{lemma}

{\it\bfseries Proof.}
In view of (\ref{2.17}), one has
\begin{align*}
&\Delta^{-1}{\rm div}(\rho u)_t+ u\cdot\nabla\Delta^{-1}{\rm div}(\rho u)-(2\mu+\lambda){\rm div}u+P+
\Delta^{-1}{\rm div}{\rm div}\bigg(\nabla d \odot \nabla d -\frac{1}{2}|\nabla d|^2\mathbb{I}_3\bigg)\\
= &u\cdot\nabla\Delta^{-1}{\rm div}(\rho u)-\Delta^{-1}{\rm div}{\rm div}(\rho u \otimes u)=[u, \mathcal{R}\otimes\mathcal{R}](\rho u),
\end{align*}
where $\mathcal{R}$ is the Riesz transform on $\mathbb{R}^3$.
To obtain the estimates of $\|\rho\|_\infty$, we adapt
the arguments by \cite{L2020}. Exactly in the same way as in Proposition 2.6 of \cite{L2020}, one can prove
Thus, we obtain
\begin{align}
\label{2.31}
\sup_{0\leq t\leq T}\|\rho\|_\infty
\leq \|\rho_0\|_\infty e^{C\overline{\rho}^{\frac{2}{3}}\sup_{0\leq t\leq T}\|\sqrt{\rho}u\|_2^{\frac{1}{3}}\|\sqrt{\rho}|u|^2\|_2^{\frac{1}{3}} +
C\overline{\rho}\int^T_0\|\nabla u\|_2\|(\nabla G, \nabla \omega, \overline{\rho}\nabla\theta)\|_2dt+ C\int^T_0\|\nabla d\|_\infty^2dt}.
\end{align}
Thanks to this and noticing that
\begin{align*}
\int^T_0\|\nabla d\|_\infty^2dt\leq C\int^T_0 \|\nabla^2 d\|_2\|\nabla^3 d\|_2dt \leq C\left(\int^T_0 \|\nabla^2 d\|_2dt\right)^\frac{1}{2}\left(\int^T_0\|\nabla^3 d\|_2dt\right)^\frac{1}{2},
\end{align*}
the conclusion follows.
\hfill$\Box$

\begin{lemma}
\label{lemma 2.8}
Let
$$
\mathcal{N}_T:=\overline{\rho}\big(\|\rho\|_3+\overline{\rho}^2(\|\sqrt{\rho}u\|_2^2+ \|\nabla d\|_2^2)\big)
\big(\|\nabla u\|_2^2+\overline{\rho}(\|\sqrt{\rho}E\|_2^2 + \|\nabla^2 d\|_2^2)\big)(t).
$$
Then there is a positive constant $\eta_0$ depending only on $R, c_v, \mu, \lambda$, and $\kappa$, such that
if
\begin{align*}
\eta\leq\eta_0, ~\sup_{0\leq t\leq T}\|\rho\|_\infty \leq 4\overline{\rho},~ \text{and~} \mathcal{N}_T\leq\sqrt{\eta},
\end{align*}
then it holds that
\begin{align}
&\sup_{0\leq t\leq T}\|\rho\|_{3}+ \left(\int^T_0\int \rho^3Pdxdt\right)^\frac{1}{3}
\leq C\big(\|\rho_0\|_3+\overline{\rho}^2(\|\sqrt{\rho_0}u_0\|_2^2+ \|\nabla d_0\|_2^2)\big), \label{2.32}\\[3mm]
&\overline{\rho}^2\left(\sup_{0\leq t\leq T}(\|\sqrt{\rho}u\|_{2}^2+ \|\nabla d\|_{2}^2)+ \int^T_0\|(\nabla u, d_t, \nabla^2 d)\|_{2}^2dt\right)\nonumber\\
&\qquad\qquad\qquad\qquad\qquad\leq C\big(\|\rho_0\|_3+\overline{\rho}^2(\|\sqrt{\rho_0}u_0\|_2^2+ \|\nabla d_0\|_2^2)\big), \label{2.33}
\\[3mm]
&\sup_{0\leq t\leq T}\big[\overline{\rho}(\| \nabla^2 d\|_{2}^2+ \|\nabla d\|_{4}^4+\|\sqrt{\rho}E\|_{2}^2)+ \|\nabla u\|_{2}^2\big]\nonumber\\
&+\int^T_0\left(\|(\nabla d_t, \nabla^3 d, |\nabla d||\nabla^2 d|, \nabla \theta, |u||\nabla u|)\|_{2}^2+\bigg\|\bigg(\sqrt{\rho}u_t, \frac{\nabla G}{\sqrt{\overline{\rho}}}, \frac{\nabla \omega}{\sqrt{\overline{\rho}}}\bigg)\bigg\|_2^2\right)dt\label{2.34}\\
&\qquad\qquad\qquad\qquad\leq
C\big(\overline{\rho}(\| \nabla^2 d_0\|_{2}^2+\|\sqrt{\rho_0}E_0\|_{2}^2)+\|\nabla u_0\|_2^2\big),\nonumber\\[3mm]
&\sup_{0\leq t\leq T}\|\rho\|_\infty
\leq \overline{\rho}e^{C\mathcal{N}_0^{\frac{1}{6}}+C\mathcal{N}_0^{\frac{1}{2}}},  \label{2.35}
\end{align}
where the constant $C>0$ depending only on $R, c_v, \mu, \lambda$, and $\kappa$.
\end{lemma}

{\it\bfseries Proof.}
In view of Lemma \ref{lemma 2.2} and by choosing $\eta_0<1$ small enough, we have
\begin{align}
\label{2.36}
\begin{aligned}
&\sup_{0\leq t\leq T}(\|\sqrt{\rho}u\|_{2}^2+ \|\nabla d\|_{2}^2)+ \int^T_0(\|\nabla u\|_{2}^2+ \|d_t\|_{2}^2+
\|\nabla^2 d\|_{2}^2)dt\\
&\leq C(\|\sqrt{\rho_0}u_0\|_{2}^2+ \|\nabla d_0\|_{2}^2)+ C\sup_{0\leq t\leq T}\|\rho\|_3^2\int^T_0\|\nabla\theta\|^2_2dt.
\end{aligned}
\end{align}
It follows from Lemma \ref{lemma 2.3} and the assumptions that
\begin{align}
\label{2.37}
\begin{aligned}
\overline{\rho}\sup_{0\leq t\leq T}(\| \nabla^2 d\|_{2}^2+ \|\nabla d\|_{4}^4)+ \overline{\rho}\int^T_0(\|\nabla d_t\|_{2}^2+ \|\nabla^3 d\|_{2}^2+
\||\nabla d||\nabla^2 d|\|_{2}^2)dt\\
\leq C\overline{\rho}(\| \nabla^2 d_0\|_{2}^2+ \|\nabla d_0\|_{4}^4)+C\eta^{\frac{1}{2}}\int^T_0\|\nabla u\|^6_2dt+C\eta^{\frac{3}{2}}\int^T_0\|\nabla^3 d\|^2_2dt.
\end{aligned}
\end{align}
With the help of (\ref{2.36}) and since $\overline{\rho}= \|\rho_0\|_\infty+1$, one deduces from the assumption that
\begin{align*}
\int^T_0\|\nabla u\|^6_2dt\leq& \sup_{0\leq t\leq T}\|\nabla u\|_2^4\int^T_0\|\nabla u\|^2_2dt\\
\leq& C\sup_{0\leq t\leq T}\|\nabla u\|_2^4\bigg(\sup_{0\leq t\leq T}(\|\sqrt{\rho}u\|_{2}^2+ \|\nabla d\|_{2}^2)+ \sup_{0\leq t\leq T}\|\rho\|_3^2\int^T_0\|\nabla\theta\|^2_2dt\bigg)\\
\leq& C\eta^{\frac{1}{2}}\sup_{0\leq t\leq T}\|\nabla u\|_2^2+
C\eta \int^T_0\|\nabla\theta\|^2_2dt,
\end{align*}
which together with (\ref{2.37}) and by choosing $\eta_0$ small enough, implies that
\begin{align}
\label{2.38}
\begin{aligned}
&\overline{\rho}\sup_{0\leq t\leq T}(\| \nabla^2 d\|_{2}^2+ \|\nabla d\|_{4}^4)+ \overline{\rho}\int^T_0(\|\nabla d_t\|_{2}^2+ \|\nabla^3 d\|_{2}^2+
\||\nabla d||\nabla^2 d|\|_{2}^2)dt\\
\leq &C\overline{\rho}(\| \nabla^2 d_0\|_{2}^2+ \|\nabla d_0\|_{4}^4)+C\eta^{\frac{1}{2}}\sup_{0\leq t\leq T}\|\nabla u\|_2^2
+C\eta\int^T_0\|\nabla\theta\|^2_2dt.
\end{aligned}
\end{align}
Next, applying Lemma \ref{lemma 2.4}, using the assumptions and (\ref{2.36}), we obtain
\begin{align*}
&\sup_{0\leq t\leq T}\|\sqrt{\rho}E\|_{2}^2+ \int^T_0(\|\nabla \theta\|_{2}^2+ \||u||\nabla u|\|_{2}^2)dt\\
\leq& C\|\sqrt{\rho_0}E_0\|_{2}^2+C\eta^\frac{1}{4} \int^T_0(\|\nabla \theta\|^2_2+ \||u||\nabla u|\|_{2}^2)dt+C\eta\int^T_0\|\nabla^3 d\|^2_2dt+ C\eta^\frac{1}{2}
\int^T_0\||\nabla d||\nabla^2 d|\|_{2}^2dt\\
&
+ C\sup_{0\leq t\leq T}\|\nabla^2 d\|_2^4(\|\sqrt{\rho}u\|_{2}^2+ \|\nabla d\|_{2}^2)+ C\sup_{0\leq t\leq T}(\|\nabla^2 d\|_2^4\|\rho\|_3^2)\int^T_0\|\nabla\theta\|^2_2dt\\
\leq& C\|\sqrt{\rho_0}E_0\|_{2}^2+C\eta^\frac{1}{4} \int^T_0(\|\nabla \theta\|^2_2+ \||u||\nabla u|\|_{2}^2)dt+C\eta\int^T_0\|\nabla^3 d\|^2_2dt+ C\eta^\frac{1}{2}
\int^T_0\||\nabla d||\nabla^2 d|\|_{2}^2dt\\
&
+ C\eta^{\frac{1}{2}}\sup_{0\leq t\leq T}\|\nabla^2 d\|_2^2+ C\eta\int^T_0\|\nabla\theta\|^2_2dt.
\end{align*}
This, combined with the fact $\eta_0$ is small enough, implies that
\begin{align}
\label{2.39}
\begin{aligned}
&\sup_{0\leq t\leq T}\|\sqrt{\rho}E\|_{2}^2+ \int^T_0(\|\nabla \theta\|_{2}^2+ \||u||\nabla u|\|_{2}^2)dt\\
\leq& C\|\sqrt{\rho_0}E_0\|_{2}^2+C\eta\int^T_0\|\nabla^3 d\|^2_2dt+ C\eta^\frac{1}{2}
\int^T_0\||\nabla d||\nabla^2 d|\|_{2}^2dt+ C\eta^{\frac{1}{2}}\sup_{0\leq t\leq T}\|\nabla^2 d\|_2^2.
\end{aligned}
\end{align}
Then, using the assumptions and Sobolev inequality, it follows from Lemma \ref{lemma 2.6} that
\begin{align}
\label{2.40}
\begin{aligned}
&\sup_{0\leq t\leq T}\|\nabla u\|_{2}^2+ \int^T_0 \|(\sqrt{\rho}u_t, \frac{\nabla G}{\sqrt{\overline{\rho}}}, \frac{\nabla \omega}{\sqrt{\overline{\rho}}})\|_2^2dt\\
\leq& C\|\nabla u_0\|_2^2+ C\overline{\rho}\sup_{0\leq t\leq T}\|\sqrt{\rho}E\|_2^2
+C\overline{\rho}^3\int^t_0\|\nabla u\|_2^4(\|\nabla u\|_2^2+\overline{\rho}\|\sqrt{\rho}E\|_2^2 )dt\\
&+C\int^T_0 (\overline{\rho} + \overline{\rho}^2\|\rho\|_3^\frac{1}{2}\|\sqrt{\rho}\theta\|_2)(\|\nabla\theta\|_2^2+\||u||\nabla u|\|_2^2)dt+ \eta\overline{\rho}\int^T_0\|\nabla\theta\|_2^2dt\\
&+C\|\nabla d_0\|_4^4+ C\|\nabla d\|_3^2\|\nabla^2 d\|_2^2 +
C\eta\overline{\rho}\int^T_0(\|\nabla d_t\|_2^2+\|\nabla^3 d\|_2^2)dt\\
&+C\eta^{\frac{1}{2}}\overline{\rho}\int^T_0\|\nabla (|\nabla d|^2) \|_2^2dt,
\end{aligned}
\end{align}
where we choose $\varepsilon_1 \leq \eta$ small enough.
By (\ref{2.36}) and (\ref{2.39}), we get
\begin{align}
\label{2.41}
\begin{aligned}
&\overline{\rho}^3\int^t_0\|\nabla u\|_2^4(\|\nabla u\|_2^2+\overline{\rho}\|\sqrt{\rho}E\|_2^2 )dt\\
\leq& C\overline{\rho}^3\sup_{0\leq t\leq T}(\|\nabla u\|_2^2+ \overline{\rho}\|\sqrt{\rho}E\|_2^2)\sup_{0\leq t\leq T}\|\nabla u\|_2^2\\
&\cdot\left((\sup_{0\leq t\leq T}(\|\sqrt{\rho} u\|_2^2+ \|\nabla d\|_2^2)+ \sup_{0\leq t\leq T}\|\rho\|_3^2\int^T_0\|\nabla\theta\|_2^2dt)
\right)\\
\leq& C\eta^{\frac{1}{2}}\sup_{0\leq t\leq T}\|\nabla u\|_2^2+ C\eta\int^T_0\|\nabla\theta\|_2^2dt
\end{aligned}
\end{align}
and
\begin{align}
\label{2.42}
\begin{aligned}
&\overline{\rho}\sup_{0\leq t\leq T}\|\sqrt{\rho}E\|_2^2+\int^T_0 (\overline{\rho} + \overline{\rho}^2\|\rho\|_3^\frac{1}{2}\|\sqrt{\rho}\theta\|_2)(\|\nabla\theta\|_2^2+\||u||\nabla u|\|_2^2)dt\\
\leq& \overline{\rho}\sup_{0\leq t\leq T}\|\sqrt{\rho}E\|_2^2+ (\overline{\rho} + \overline{\rho}^2
\sup_{0\leq t\leq T}(\|\rho\|_3^\frac{1}{2}\|\sqrt{\rho}E\|_2)\int^T_0 (\|\nabla\theta\|_2^2+\||u||\nabla u|\|_2^2)dt\\
\leq& \overline{\rho}\sup_{0\leq t\leq T}\|\sqrt{\rho}E\|_2^2+ (\overline{\rho} + \overline{\rho}
\eta^\frac{1}{4})\int^T_0 (\|\nabla\theta\|_2^2+\||u||\nabla u|\|_2^2)dt\\
\leq& C\overline{\rho}\|\sqrt{\rho_0}E_0\|_{2}^2+C\overline{\rho}\eta\int^T_0\|\nabla^3 d\|^2_2dt+ C\overline{\rho}\eta^\frac{1}{2}
\int^T_0\||\nabla d||\nabla^2 d|\|_{2}^2dt+ C\overline{\rho}\eta^{\frac{1}{2}}\sup_{0\leq t\leq T}\|\nabla^2 d\|_2^2.
\end{aligned}
\end{align}
Substituting (\ref{2.41}) and (\ref{2.42}) into (\ref{2.40}) and using $\eta_0$ is small enough, one obtains
\begin{align}
\label{2.43}
\begin{aligned}
&\sup_{0\leq t\leq T}\|\nabla u\|_{2}^2+ \int^T_0 \|(\sqrt{\rho}u_t, \frac{\nabla G}{\sqrt{\overline{\rho}}}, \frac{\nabla \omega}{\sqrt{\overline{\rho}}})\|_2^2dt\\
\leq& C(\|\nabla u_0\|_2^2+\overline{\rho}\|\sqrt{\rho_0}E_0\|_{2}^2+\|\nabla d_0\|_4^4)
+ C\eta\overline{\rho}\int^T_0\|\nabla\theta\|_2^2dt\\
&+ C\overline{\rho}\eta^\frac{1}{2}
\int^T_0\||\nabla d||\nabla^2 d|\|_{2}^2dt+ C\overline{\rho}\eta^{\frac{1}{2}}\sup_{0\leq t\leq T}\|\nabla^2 d\|_2^2 \\
&+C\eta\overline{\rho}\int^T_0\|\nabla d_t\|_2^2dt
+C\eta\overline{\rho}\int^T_0\|\nabla^3 d\|_2^2dt.
\end{aligned}
\end{align}

The combination of (\ref{2.38}), (\ref{2.39}) and (\ref{2.43}) yields that
\begin{align*}
&\sup_{0\leq t\leq T}\big(\overline{\rho}(\| \nabla^2 d\|_{2}^2+ \|\nabla d\|_{4}^4+\|\sqrt{\rho}E\|_{2}^2)+ \|\nabla u\|_{2}^2\big)+\int^T_0\bigg\|\bigg(\sqrt{\rho}u_t, \frac{\nabla G}{\sqrt{\overline{\rho}}}, \frac{\nabla \omega}{\sqrt{\overline{\rho}}}\bigg)\bigg\|_2^2dt\\
&+\overline{\rho}\int^T_0\left(\|\nabla d_t\|_{2}^2+ \|\nabla^3 d\|_{2}^2+\||\nabla d||\nabla^2 d|\|_{2}^2+
\|\nabla \theta\|_{2}^2+ \||u||\nabla u|\|_{2}^2\right) dt\\
\leq&
C\big(\overline{\rho}(\| \nabla^2 d_0\|_{2}^2+\|\sqrt{\rho_0}E_0\|_{2}^2)+\|\nabla u_0\|_2^2\big)\\
&+ C\overline{\rho}\sup_{0\leq t\leq T}\|\nabla d\|_{3}^2\|\nabla^2 d\|_{2}^2+C\eta^{\frac{1}{2}}\sup_{0\leq t\leq T}\|\nabla u\|_2^2
+C\eta\overline{\rho}\int^T_0\|\nabla\theta\|^2_2dt\\
&+ C\overline{\rho}\eta^\frac{1}{2}
\int^T_0\||\nabla d||\nabla^2 d|\|_{2}^2dt+ C\overline{\rho}\eta^{\frac{1}{2}}\sup_{0\leq t\leq T}\|\nabla^2 d\|_2^2 +
C\eta\overline{\rho}\int^T_0\|\nabla d_t\|_2^2dt
+C\eta\overline{\rho}\int^T_0\|\nabla^3 d\|_2^2dt\\
\leq&
C\big(\overline{\rho}(\| \nabla^2 d_0\|_{2}^2+\|\sqrt{\rho_0}E_0\|_{2}^2)+\|\nabla u_0\|_2^2\big)+C\eta^{\frac{1}{2}}\sup_{0\leq t\leq T}\|\nabla u\|_2^2
+C\eta\overline{\rho}\int^T_0\|\nabla\theta\|^2_2dt\\
&+ C\overline{\rho}\eta^\frac{1}{2}
\int^T_0\||\nabla d||\nabla^2 d|\|_{2}^2dt+ C\overline{\rho}(\eta^{\frac{1}{2}}+\eta)\sup_{0\leq t\leq T}\|\nabla^2 d\|_2^2 \\
&+
C\eta\overline{\rho}\int^T_0\|\nabla d_t\|_2^2dt
+C\eta\overline{\rho}\int^T_0\|\nabla^3 d\|_2^2dt,
\end{align*}
from which, choosing $\eta_0$ small enough, one gets (\ref{2.34}) and
\begin{align}
\label{2.44}
\overline{\rho}\int^T_0\|\nabla\theta\|^2_2dt\leq C\big(\overline{\rho}(\| \nabla^2 d_0\|_{2}^2+\|\sqrt{\rho_0}E_0\|_{2}^2)+\|\nabla u_0\|_2^2\big).
\end{align}
Recalling (\ref{2.36}) and the assumptions, and using (\ref{2.44}), we have
\begin{align}
\label{2.45}
\begin{aligned}
&\sup_{0\leq t\leq T}(\|\sqrt{\rho}u\|_{2}^2+ \|\nabla d\|_{2}^2)+ \int^T_0(\|\nabla u\|_{2}^2+ \|d_t\|_{2}^2+
\|\nabla^2 d\|_{2}^2)dt\\
&\leq C(\|\sqrt{\rho_0}u_0\|_{2}^2+ \|\nabla d_0\|_{2}^2)+ C\frac{1}{\overline{\rho}}\sup_{0\leq t\leq T}\|\rho\|_3^2\big(\overline{\rho}(\| \nabla^2 d_0\|_{2}^2+\|\sqrt{\rho_0}E_0\|_{2}^2)+\|\nabla u_0\|_2^2\big)\\
&\leq C(\|\sqrt{\rho_0}u_0\|_{2}^2+ \|\nabla d_0\|_{2}^2)+ C\eta^{\frac{1}{2}}\frac{1}{\overline{\rho}^2}\sup_{0\leq t\leq T}\|\rho\|_3.
\end{aligned}
\end{align}

It follows from Lemma \ref{lemma 2.5}, (\ref{2.45}), the Young inequality and the assumptions that
\begin{align*}
&\sup_{0\leq t\leq T}\|\rho\|_{3}^3+ \int^T_0\int \rho^3Pdxdt\\
\leq& C\|\rho_0\|_3^3+ C\sup_{0\leq t\leq T}(\|\rho\|_\infty^{\frac{2}{3}}\|\sqrt{\rho}u\|_2^{\frac{1}{3}}\|\sqrt{\rho}E\|_2^{\frac{1}{3}}
\|\rho\|_3^3)+C\overline{\rho}^2\sup_{0\leq t\leq T}\|\rho\|_3^2\int^T_0(\|\nabla u\|^2_2+\|\nabla^2 d\|_2^2)dt
\\
\leq& C\|\rho_0\|_3^3+ C\eta^{\frac{1}{12}}\sup_{0\leq t\leq T}\|\rho\|_3^3 +C\overline{\rho}^2\sup_{0\leq t\leq T}\|\rho\|_3^2\left(\|\sqrt{\rho_0}u_0\|_{2}^2+ \|\nabla d_0\|_{2}^2+ \eta^{\frac{1}{2}}\frac{1}{\overline{\rho}^2}\sup_{0\leq t\leq T}\|\rho\|_3\right)\\
\leq& C\|\rho_0\|_3^3+ C\big(\eta^{\frac{1}{12}}+\eta^{\frac{1}{2}} \big)\sup_{0\leq t\leq T}\|\rho\|_3^3 +C\overline{\rho}^2\sup_{0\leq t\leq T}\|\rho\|_3^2(\|\sqrt{\rho_0}u_0\|_{2}^2+ \|\nabla d_0\|_{2}^2)\\
\leq& C\|\rho_0\|_3^3+ C\big(\eta^{\frac{1}{12}}+\eta^{\frac{1}{2}} +\frac{1}{4}\big)\sup_{0\leq t\leq T}\|\rho\|_3^3 +C\overline{\rho}^6(\|\sqrt{\rho_0}u_0\|_{2}^2+ \|\nabla d_0\|_{2}^2)^3,
\end{align*}
which implies (\ref{2.32}) by choosing $\eta_0$ sufficiently small.

Now, substituting (\ref{2.32}) into (\ref{2.45}) yields that
\begin{align*}
&\overline{\rho}^2\left(\sup_{0\leq t\leq T}(\|\sqrt{\rho}u\|_{2}^2+ \|\nabla d\|_{2}^2)+ \int^T_0(\|\nabla u\|_{2}^2+ \|d_t\|_{2}^2+
\|\nabla^2 d\|_{2}^2)dt\right)\\
&\leq C\overline{\rho}^2(\|\sqrt{\rho_0}u_0\|_{2}^2+ \|\nabla d_0\|_{2}^2)+ C\eta^{\frac{1}{2}}\sup_{0\leq t\leq T}\|\rho\|_3\\
&\leq C\big(\|\rho_0\|_3+ \overline{\rho}^2(\|\sqrt{\rho_0}u_0\|_{2}^2+ \|\nabla d_0\|_{2}^2)\big),
\end{align*}
which gives (\ref{2.33}).

Finally, (\ref{2.35}) follows immediately from Lemma \ref{lemma 2.7}, (\ref{2.33}) and (\ref{2.34}), and the proof is complete.
\hfill$\Box$

\begin{lemma}
\label{lemma 2.9}
Let $\eta_0, \mathcal{N}_T$, and $\mathcal{N}_0$ be as in Lemma \ref{lemma 2.8}. Then, there exists a number $\varepsilon_0\in (0, \eta_0)$
such that if
$$
\sup_{0\leq t\leq T}\|\rho\|_{\infty} \leq 4\overline{\rho}, ~\mathcal{N}_T\leq\sqrt{\varepsilon_0} \text{~~and~~}
\mathcal{N}_0\leq \varepsilon_0,
$$
then
$$
\sup_{0\leq t\leq T}\|\rho\|_{\infty} \leq 2\overline{\rho} \text{~~and~~} \mathcal{N}_T\leq\frac{\sqrt{\varepsilon_0}}{2},
$$
where $\varepsilon_0$ depending only on $R, c_v, \mu, \lambda$, and $\kappa$.
\end{lemma}

{\it\bfseries Proof.}
If $\varepsilon_0\leq \eta_0$ is sufficiently small, all the conditions in Lemma \ref{lemma 2.8} hold.
Therefore, we obtain
\begin{align*}
\mathcal{N}_T \leq C\overline{\rho}\big(\|\rho_0\|_3+\overline{\rho}^2(\|\sqrt{\rho_0}u_0\|_2^2+ \|\nabla d_0\|_2^2)\big)
\big(\|\nabla u_0\|_2^2+\overline{\rho}(\|\sqrt{\rho_0}E_0\|_2^2 + \|\nabla^2 d_0\|_2^2)\big)
\leq C\varepsilon_0 \leq \frac{\sqrt{\varepsilon_0}}{2}.
\end{align*}
At the same time,
\begin{align*}
\sup_{0\leq t\leq T}\|\rho\|_\infty
\leq \overline{\rho}e^{C\mathcal{N}_0^{\frac{1}{6}}+C\mathcal{N}_0^{\frac{1}{2}}}
\leq\overline{\rho}e^{C\varepsilon_0^{\frac{1}{6}}+C\varepsilon_0^{\frac{1}{2}}}\leq 2\overline{\rho}.
\end{align*}
We complete the proof of the lemma.
\hfill$\Box$

Thus, based on Lemma \ref{lemma 2.1} and Lemma \ref{lemma 2.9}, together with the
standard continuity argument, one can deduce the following
proposition:

\begin{proposition}
\label{proposition 2.10}
Assume $\mathcal{N}_0\leq \varepsilon_0$ with $\varepsilon_0$ defined as in Lemma \ref{lemma 2.9}. Then,
it holds that
\begin{align*}
\mathcal{N}_T\leq \frac{\varepsilon_0}{2} \text{~~and~} \sup_{0\leq t\leq T}\|\rho\|_\infty \leq2\overline{\rho}.
\end{align*}
Moreover, there is
a constant $C>0$ such that the following estimates hold:
\begin{align*}
&\sup_{0\leq t\leq T}\big(\|(\sqrt{\rho}E, \sqrt{\rho}u, \nabla u, d_t, \nabla d, \nabla^2 d)\|_2^2+
\|\rho\|_3+ \|\rho\|_\infty+ \|\nabla d\|_4^4+ \|\nabla d\|_3^3 \big)\leq C,\\
&\int^T_0 \big(\|(\nabla\theta, |u|\nabla u, \sqrt{\rho}u_t, \nabla u, \nabla G, \nabla\omega, d_t, \nabla d_t, \nabla^2d,
 \nabla^3d, |\nabla d||\nabla^2d|)\|_2^2+ \|\nabla u\|_6^2\big)dt\\
 & + \int^T_0\int \rho^3Pdxdt\leq C,
\end{align*}
where $C$ depends only on $R, c_v, \mu, \lambda$, $\kappa, \overline{\rho}, \|\rho_0\|_3, \|\sqrt{\rho_0}u_0\|_2, \|\sqrt{\rho_0}E_0\|_2,
\|\nabla u_0\|_2, \|\nabla d_0\|_2, \|\nabla^2 d_0\|_2$, and $\|\nabla d_0\|_3$.
\end{proposition}

\section{Time dependent higher order estimates}

Taking advantage of Proposition \ref{proposition 2.10}, we can obtain
the higher order estimates which is sufficient to ensure the existence of global strong solutions.
Precisely, the following estimate is proved in this section
$$
\sup_{0\leq t\leq T}(\|\nabla\theta\|_{H^1}^2+\|(\nabla^2u, \sqrt\rho\dot u, \sqrt\rho\dot\theta, \nabla^3d, \nabla d_t)\|_2^2+\|\rho\|_{H^1\cap W^{1,q}}) \leq C.
$$
This a priori estimates can be established by modifying the methods of \cite{HL2013, HL2018, XZ2012}
for the compressible Navier-Stokes equations and magnetohydrodynamic equations.

In the rest of this section, we always assume that $(\rho, u, \theta, d)$ is a strong solution to system (\ref{1.1})-(\ref{1.5}) in $\mathbb{R}^3\times (0, T)$, for a positive time $T$.

\begin{lemma}
\label{lemma 3.1}
Assume $\mathcal{N}_0\leq \varepsilon_0$. It holds that
\begin{align*}
\sup_{0\leq t\leq T}(\|\nabla \theta\|_2^2+\|\sqrt{\rho}\dot{u}\|_2^2+\|\nabla d_t\|_2^2)
+\int^T_0 \|(\sqrt{\rho}\dot{\theta}, \nabla\dot{u}, d_{tt}, \Delta d_t)\|_2^2dt\leq C_T,
\end{align*}
where $C_T$ depending only on $R, c_v, \mu, \lambda, \kappa, \Phi_0,$ and $T$.
\end{lemma}

{\it\bfseries Proof.}
Applying $\dot{u}_j\big( \partial_t + {\rm div}(u\cdot)\big)$ to $(\ref{1.2})^{j}$ and integrating over $\mathbb{R}^3$, it follows
\begin{align}
\label{3.1}
\begin{aligned}
\frac{1}{2}\frac{d}{dt} \|\sqrt{\rho}\dot{u}\|_2^2=&-\int \dot{u}_j\big( \partial_jP_t + {\rm div}(u \partial_jP)\big)dx
+\mu \int \dot{u}_j\big( \partial_t\Delta u_j + {\rm div}(u \Delta u_j)\big)dx\\
&+(\mu+\lambda) \int \dot{u}_j\big( \partial_j{\rm div}u_t + {\rm div}(u \partial_j{\rm div}u)\big)dx\\
&-\int \partial_i(M_{i,j}(d))_t \cdot \dot{u}_jdx -\int \partial_k\big(u_k\partial_i(M_{i,j}(d))\big)\dot{u}_j dx
=:\sum_{i=1}^5J_i,
\end{aligned}
\end{align}
where $M_{i,j}(d) = \partial_i d\cdot\partial_j d-\frac{1}{2}|\nabla d|^2\delta_{i,j}$.
It follows from the H\"{o}lder, Young and Sobolev inequalities that
\begin{align*}
J_1=&-\int \dot{u}_j\big[ \partial_jP_t + \partial_j {\rm div}(u P)- {\rm div}(\partial_juP)\big]dx\\
=&\int {\rm div}\dot{u}(P_t+ {\rm div}(uP))dx- \int \nabla\dot{u}_j\cdot \partial_juPdx\\
=& R\int{\rm div}\dot{u} \rho\dot{\theta} dx - R\int \nabla\dot{u}_j\cdot \partial_ju \rho\theta dx\\
\leq &\frac{\mu}{8}\|\nabla \dot{u}\|_2^2+ C\|\rho\dot{\theta}\|_2^2
+C\int \rho^2\theta^2|\nabla u|^2dx\\
\leq &\frac{\mu}{8}\|\nabla \dot{u}\|_2^2+ C\|\rho\dot{\theta}\|_2^2
+C\|\rho\theta\|_2^{\frac{1}{2}}\|\theta\|_6^{\frac{3}{2}}\|\nabla u\|_4^2\\
\leq &\frac{\mu}{8}\|\nabla \dot{u}\|_2^2+ C\big(1+\|\rho\dot{\theta}\|_2^2
+\|\nabla\theta\|_2^4+ \|\nabla u\|_4^{4}\big),
\end{align*}
where Proposition \ref{proposition 2.10} was used.
By virtue of integration by parts, we compute
\begin{align*}
J_2=&-\mu\int\big(\partial_i\dot{u}_j(\partial_iu_j)_t+\Delta u_ju\cdot\nabla \dot{u}_j\big)dx\\
=&-\mu\int\big(|\nabla\dot{u}|^2-\partial_i\dot{u}_ju_k\partial_k\partial_iu_j
-\partial_i\dot{u}_j\partial_iu_k\partial_ku_j+\Delta u_ju\cdot\nabla\dot{u}_j
\big)dx
\\
=&-\mu\int\big(|\nabla\dot{u}|^2+\partial_i\dot{u}_j \partial_iu_j{\rm div}u
-\partial_i\dot{u}_j\partial_iu_k\partial_ku_j-\partial_i u_j\partial_iu_k\partial_k \dot{u}_j
\big)dx\\
\leq &-\frac{7\mu}{8}\|\nabla \dot{u}\|_2^2+ C\|\nabla u\|_4^4.
\end{align*}
In the same way, one gets
\begin{align*}
J_3
\leq -\frac{7(\mu+\lambda)}{8}\|{\rm div} \dot{u}\|_2^2+ C\|\nabla u\|_4^4.
\end{align*}
For $J_4$ and $J_5$, by integration by parts and the H\"{o}lder, Young and Sobolev inequalities, one has for $\eta_1\in (0,1]$
\begin{align*}
J_4
\leq& \int |\nabla d||\nabla d_t||\nabla\dot{u}|dx\leq C\|\nabla d\|_6\|\nabla d_t\|_3\|\nabla\dot{u}\|_2\leq C\|\nabla^2 d\|_2\|\nabla d_t\|_2^{\frac{1}{2}}\|\nabla^2 d_t\|_2^{\frac{1}{2}}
\|\nabla\dot{u}\|_2\\
\leq& \varepsilon\|\nabla\dot{u}\|_2^2+ \eta_1\|\nabla^2 d_t\|_2^2 +C(\varepsilon,\eta_1)\|\nabla d_t\|_2^2,
\end{align*}
where $\displaystyle\sup_{0\leq t\leq T}\|\nabla^2d\|_2 \leq C$ guaranteed by \ref{proposition 2.10} was used, and
\begin{align*}
J_5
\leq& \int |u||\nabla d||\nabla^2d||\nabla\dot{u}|dx\leq C\|u\|_6\|\nabla d\|_6\|\nabla^2 d\|_6\|\nabla\dot{u}\|_2\leq
\varepsilon\|\nabla\dot{u}\|_2^2+ C(\varepsilon)\|\nabla^3 d\|_2^2,
\end{align*}
where Proposition \ref{proposition 2.10} was used.

Substituting $J_i, i=1,2, \ldots, 5$ into (\ref{3.1}), one obtains after choosing $\varepsilon$ suitably small that
\begin{align}
\label{3.2}
\begin{aligned}
&\frac{d}{dt} \|\sqrt{\rho}\dot{u}\|_2^2+ \mu\|\nabla\dot{u}\|_2^2\\
\leq&C\eta_1\|\nabla^2 d_t\|_2^2 + C\big(1+\|\rho\dot{\theta}\|_2^2
+\|\nabla\theta\|_2^4+ \|\nabla u\|_4^{4}+ \|\nabla^3 d\|_2^2\big)+
C(\eta_1)\|\nabla d_t\|_2^2.
\end{aligned}
\end{align}

Next, multiplying (\ref{1.3}) by $\dot{\theta}$ and integrating the resultant over $\mathbb{R}^3$ yield
\begin{align}
\label{3.3}
\begin{aligned}
\frac{\kappa}{2}\frac{d}{dt}\|\nabla\theta\|_2^2 + c_v\|\sqrt{\rho}\dot{\theta}\|_2^2
=&
-\kappa\int \nabla\theta\cdot\nabla(u\cdot\nabla\theta)dx+ \lambda\int |{\rm div}u|^2\dot{\theta}dx+ \frac{\mu}{2}\int |\nabla u+(\nabla u)^t|^2\dot{\theta}dx\\
&-R\int \rho\theta{\rm div}u\dot{\theta} dx+ \int |\Delta d+|\nabla d|^2d|^2\dot{\theta}dx:=\sum_{i=1}^5K_i.
\end{aligned}
\end{align}
It follows from elliptic estimates, Proposition \ref{proposition 2.10}, and Gagliardo-Nirenberg and Young inequalities that
\begin{align*}
\|\nabla \theta\|_{H^1}^2\leq& C+ C\|\sqrt{\rho}\dot{\theta}\|_2^2+ C\|\nabla \theta\|_2^2+ C\int \rho^2\theta^2|\nabla u|^2dx+C\|\nabla u\|_4^4+ C\|\Delta d+|\nabla d|^2d\|_4^4\\
\leq& C+ C\|\sqrt{\rho}\dot{\theta}\|_2^2+ C\|\nabla \theta\|_2^2+ C\|\nabla u\|_2^2\|\theta\|_\infty^2+C\|\nabla u\|_4^4+ C\|\Delta d+|\nabla d|^2d\|_4^4\\
\leq& \frac{1}{2}\|\nabla\theta\|_{H^1}^2+C\big(1+\|\sqrt{\rho}\dot{\theta}\|_2^2+ \|\nabla \theta\|_2^2+\|\nabla u\|_4^4+ \|\Delta d+|\nabla d|^2d\|_4^4\big),
\end{align*}
which implies
\begin{align}
\label{3.4}
\begin{aligned}
\|\nabla \theta\|_{H^1}^2\leq C\big(1+ \|\sqrt{\rho}\dot{\theta}\|_2^2+ \|\nabla \theta\|_2^2+\|\nabla u\|_4^4+ \|\Delta d+|\nabla d|^2d\|_4^4\big).
\end{aligned}
\end{align}
Moveover, by (\ref{1.4}), the H\"{o}lder, Young and Sobolev inequalities and Proposition \ref{proposition 2.10}, one can get by the elliptic estimates that
\begin{align*}
\|\nabla^3 d\|_2 \leq& C\|\nabla d_t\|_2+ \|\nabla u\cdot\nabla d\|_2+ C\|u\cdot\nabla^2d\|_2+C\||\nabla d|^3\|_2+ C\||\nabla d|\nabla^2d\|_2\\
 \leq& C\|\nabla d_t\|_2+C\|\nabla u\|_4\|\nabla d\|_4 + C\|u\|_6\|\nabla^2 d\|_3
 +C\|\nabla d\|_6^3+ C\|\nabla d\|_6\|\nabla^2d\|_3\\
 \leq& C\|\nabla d_t\|_2+C\|\nabla u\|_4 + C\|\nabla^2 d\|_2^{\frac{1}{2}}\|\nabla^3 d\|_2^{\frac{1}{2}}
 +C\\
  \leq& \frac{1}{2}\|\nabla^3d\|_2+C\big(1+\|\nabla d_t\|_2+\|\nabla u\|_4\big),
\end{align*}
which gives
\begin{align}
\label{3.5}
\|\nabla d\|_{H^2}\leq C\big(1+\|\nabla d_t\|_2+\|\nabla u\|_4\big).
\end{align}
According to (\ref{3.5}), one gets by the Sobolev embedding inequality and Proposition \ref{proposition 2.10} that
\begin{align}
\label{3.6}
\begin{aligned}
\|\Delta d+|\nabla d|^2d\|_4^4\leq C(\|\Delta d\|_{H^1}^4+\|\nabla d\|_\infty^4\|\nabla d\|_4^4)\leq
 C\|\nabla d\|_{H^2}^4
\leq C\big(1+\|\nabla d_t\|_{2}^4+ \|\nabla u\|_{4}^4\big).
\end{aligned}
\end{align}
Thus, (\ref{3.4}), (\ref{3.6}), the Sobolev and Young inequalities yield
\begin{align*}
K_1 =&-\kappa\int \nabla\theta\cdot\nabla(u\cdot\nabla\theta)dx\leq C\int |\nabla\theta|\big(|u||\nabla^2\theta|+|\nabla u||\nabla\theta|\big) dx\\
\leq& C\big(\|\nabla\theta\|_3\|u\|_6\|\nabla^2\theta\|_2+ \|\nabla u\|_2\|\nabla\theta\|_6\|\nabla\theta\|_3\big)\\
\leq& C\|\nabla u\|_2
\|\nabla\theta\|_2^{\frac{1}{2}}\|\nabla^2\theta\|_2^{\frac{3}{2}}
\leq \varepsilon\|\nabla\theta\|_2^2 + C(\varepsilon)\|\nabla\theta\|_2^2\\
\leq& C\varepsilon\|\sqrt{\rho}\dot{\theta}\|_2^2 +
C(\varepsilon)\big(1+\|\nabla\theta\|_2^2+\|\nabla u\|_4^4+ \|\nabla d_t\|_{2}^4\big).
\end{align*}
By integration by parts, it follows from the H\"{o}lder, Young and Sobolev inequalities and Proposition \ref{proposition 2.10} that
\begin{align*}
K_2=&\lambda\int({\rm div}u)^2\theta_tdx+\lambda\int({\rm div}u)^2u\cdot\nabla\theta dx\\
=&\lambda\left(\int({\rm div}u)^2\theta dx\right)_t-2\lambda\int \theta{\rm div}u
{\rm div}(\dot{u}-u\cdot\nabla u)dx+ \lambda\int ({\rm div}u)^2u\cdot\nabla\theta dx\\
=&\lambda\left(\int({\rm div}u)^2\theta dx\right)_t-2\lambda\int \theta{\rm div}u
{\rm div}\dot{u}dx+ 2\lambda\int\theta{\rm div}u\partial_iu_j\partial_ju_i dx
+\lambda\int u\cdot\nabla\big(\theta({\rm div}u)^2\big)dx\\
\leq&\lambda\left(\int({\rm div}u)^2\theta dx\right)_t
+C\|\theta\|_6\|\nabla u\|_2^{\frac{1}{3}}\|\nabla u\|_4^{\frac{2}{3}}(\|\nabla \dot{u}\|_2+\|\nabla u\|_4^2)\\
\leq&\lambda\left(\int({\rm div}u)^2\theta dx\right)_t
+\eta_1\|\nabla \dot{u}\|_2^2+C(\eta_1)\big(1+\|\nabla u\|_4^4+ \|\nabla\theta\|_2^4\big).
\end{align*}
Similarly, one has
\begin{align*}
K_3
\leq \frac{\mu}{2}\left(\int|\nabla u+(\nabla u)^t|^2\theta dx\right)_t
+\eta_1\|\nabla \dot{u}\|_2^2+C(\eta_1)\big(1+\|\nabla u\|_4^4+ \|\nabla\theta\|_2^4\big).
\end{align*}
Using Proposition \ref{proposition 2.10} again, we get
\begin{align*}
K_4
\leq
C\| \sqrt{\rho} \dot{\theta}\|_2 \| \sqrt{\rho} \theta\|_2^\frac{1}{4}   \|\theta\|_6^{\frac{3}{4}}\|\nabla u\|_4
\leq \varepsilon\|\sqrt{\rho}\dot{\theta}\|_2^2+C(\varepsilon)\big(1+ \|\nabla\theta\|_2^4+ \|\nabla u\|_4^4\big).
\end{align*}
At last, for $K_5$, noticing that $|\Delta d+|\nabla d|^2d|^2=|\Delta d|^2-|\nabla d|^4$ and $\Delta d\cdot d=-|\nabla d|^2$ guaranteed by $|d|=1$, it follows from the H\"{o}lder, Young and Sobolev inequalities that
\begin{align*}
K_5=&\int |\Delta d+|\nabla d|^2d|^2 \theta_tdx+\int |\Delta d+|\nabla d|^2d|^2 u\cdot\nabla\theta dx\\
=& \left(\int |\Delta d+|\nabla d|^2d|^2 \theta dx\right)_t -2\int\big(\Delta d\cdot\Delta d_t- 4|\nabla d|^2\nabla d:\nabla d_t\big)\theta dx\\
&+\int |\Delta d- (\Delta d\cdot d)d|^2 u\cdot\nabla\theta dx\\
=& \left(\int |\Delta d+|\nabla d|^2d|^2 \theta dx\right)_t -2\int\big(\Delta d\cdot\Delta d_t + 4 \Delta d\cdot d\nabla d:\nabla d_t\big)\theta dx\\
&+\int |\Delta d- (\Delta d\cdot d)d|^2 u\cdot\nabla\theta dx\\
\leq & \left(\int |\Delta d+|\nabla d|^2d|^2 \theta dx\right)_t + C\|\theta\|_6\|\Delta d_t\|_2\|\Delta d\|_3+ C\|\Delta d\|_2 \|\nabla d\|_6 \|\nabla d_t\|_6 \|\theta\|_6
\\
&+ C\|\nabla\theta\|_2\|u\|_6\|\nabla^2d\|_6^2\\
\leq & \left(\int |\Delta d+|\nabla d|^2d|^2 \theta dx\right)_t +C\|\nabla\theta\|_2\|\Delta d_t\|_2\|\nabla^2 d\|_2^{\frac{1}{2}}\|\nabla^2 d\|_2^{\frac{1}{2}}+ C\|\Delta d_t\|_2 \|\nabla\theta\|_2\\
&+C \|\nabla\theta\|_2\|\nabla u\|_2\|\nabla^3 d\|_2^2\\
\leq & \left(\int |\Delta d+|\nabla d|^2d|^2 \theta dx\right)_t +\eta_1\|\Delta d_t\|_2^2 +C(\eta_1)\big(1+\|\nabla\theta\|_2^2)(1+\|\nabla^3d\|_2^2\big).
\end{align*}
Now, substituting the estimates for $K_i, i=1,2,\ldots, 5$ into (\ref{3.3}) and then choosing $\varepsilon$ small enough,
we deduce that
\begin{align}
\label{3.7}
\begin{aligned}
\frac{d}{dt}\int \Phi dx + c_v\|\sqrt{\rho}\dot{\theta}\|_2^2\leq &
C(1+\|\nabla\theta\|_2^2)(1+\|\nabla^3d\|_2^2+\|\nabla\theta\|_2^2)+C\eta_1\|\nabla\dot{u}\|_2^2+C\eta_1\| \Delta d_t\|_2^2\\
&+
C\|\nabla u\|_4^4 + C\|\nabla d_t\|_2^4 +C,
\end{aligned}
\end{align}
where
\begin{align}
\label{3.7-1}
\Phi:=\kappa|\nabla \theta|^2 - 2\theta \big(\lambda ({\rm div}u)^2+ \frac{\mu}{2}|\nabla u+(\nabla u)^t|^2+ |\Delta d+|\nabla d|^2d|^2\big).
\end{align}

On the other hand, applying $\partial_t$ to (\ref{1.4}), we have
\begin{align*}
d_{tt}-\Delta d_t=(-u\cdot\nabla d+ |\nabla d|^2d)_t.
\end{align*}
It follows from integration by parts, the Sobolev and Young inequalities and Proposition \ref{proposition 2.10} that
\begin{align*}
\frac{d}{dt}\|\nabla d_t\|_2^2 + \|(d_{tt}, \Delta d_t)\|_2^2=&\int |d_{tt}-\Delta d_t|^2dx\\
=&\int |(-u\cdot \nabla d)_t +(|\nabla d|^2d)_t|^2dx\\
\leq & \int \big( |u_t|^2|\nabla d|^2 + |u|^2|\nabla d_t|^2+ |\nabla d|^4|d_t|^2 + |\nabla d|^2|\nabla d_t|^2 \big)dx\\
\leq & \int \big( |\dot{u}|^2|\nabla d|^2 + |u|^2|\nabla u|^2|\nabla d|^2 \big)dx+C\|\nabla d\|_6^4\|d_t\|_6^2 \\
 &+C(\|u\|_6^2+\|\nabla d\|_6^2)\|\nabla d_t\|_2\|\nabla d_t\|_6
\\
\leq & \frac{1}{2}\|\Delta d_t\|_2^2+ C\|\dot{u}\|_6^2\|\nabla d\|_3^2+ C\|u\|_6^2\|\nabla u\|_6^2\|\nabla d\|_6^2+C\|\nabla d_t\|_2^2\\
\leq & \frac{1}{2}\|\Delta d_t\|_2^2+ C\|\nabla\dot{u}\|_2^2+ C\|\nabla u\|_6^2+C\|\nabla d_t\|_2^2,
\end{align*}
which yields
\begin{align}
\label{3.8}
\frac{d}{dt}\|\nabla d_t\|_2^2 + \|d_{tt}\|_2^2+ \frac{1}{2}\|\Delta d_t\|_2^2\leq  C\|\nabla\dot{u}\|_2^2+ C\|\nabla u\|_6^2+C\|\nabla d_t\|_2^2.
\end{align}

Thus, multiplying (\ref{3.2}) and (\ref{3.8}) by $\eta_1^{\frac{1}{4}}$ and $\eta_1^{\frac{1}{2}}$, respectively, then adding the result
with (\ref{3.7}) and choosing $\eta_1$ suitably small, we finally obtain
\begin{align}
\label{3.9}
\begin{aligned}
&2\frac{d}{dt}\int \big(\Phi+\eta_1^{\frac{1}{2}}|\nabla d_t|^2 + \eta_1^{\frac{1}{4}}\rho| \dot{u}|^2\big)dx
 + c_v\|\sqrt{\rho}\dot{\theta}\|_2^2+ \eta_1^{\frac{1}{2}}\|(d_{tt}, \Delta d_t)\|_2^2+ \mu\eta_1^{\frac{1}{4}}\|\nabla\dot{u}\|_2^2\\
\leq &
C\big(1+\|\nabla\theta\|_2^2)(1+\|\nabla^3d\|_2^2+\|\nabla\theta\|_2^2\big)
+
C\big(\|\nabla u\|_4^4 + \|\nabla d_t\|_2^4 +\|\nabla d_t\|_2^4+\|\nabla u\|_6^2\big),
\end{aligned}
\end{align}
where $\Phi$ is given by (\ref{3.7-1}).

Next, one needs to show the estimate of $\|\nabla u\|_6$ in order to bound $\|\nabla u\|_4$. To this end, decompose $u=v+w$, where $v$ satisfies
\begin{align}
\label{3.10}
\mu \Delta v+(\mu+\lambda)\nabla {\rm div}v=\nabla P.
\end{align}
According to Lemma 2.3 in \cite{HL2013},
there exists a unique $v(\cdot, t)\in D^1_0\cap D^{2,2}\cap D^{2,{q}}$ satisfying (\ref{3.10}) and the following $L^p$, $p\in [2,6]$ and $L^\infty$ estimates for $t\in[0,T]$:
\begin{align}
\label{3.11}
\|\nabla v\|_p\leq C\|\rho\theta\|_p,
\end{align}
and
\begin{align}
\label{3.111}
\|\nabla v\|_\infty \leq C\big(1+\log(e+\|\nabla(\rho\theta)\|_{{q}})\|\rho\theta\|_\infty + \|\rho\theta\|_2\big), ~~{q}\in (3,6].
\end{align}
While, $w$ satisfies
\begin{align}
\label{3.12}
\mu \Delta w+(\mu+\lambda)\nabla {\rm div}w= \rho\dot{u}+ \Delta d\cdot\nabla d.
\end{align}
By the elliptic estimates, it holds that
\begin{align}
\label{3.13}
\|\nabla w\|_6 + \|\nabla^2 w\|_2\leq C\|\rho\dot{u}\|_2+C\|\Delta d\cdot\nabla d\|_2,
\end{align}
and
\begin{align}
\label{3.222}
\|\nabla^2 w\|_6 \leq C\|\rho\dot{u}\|_6+C\|\Delta d\cdot\nabla d\|_6.
\end{align}
Hence, by the H\"{o}lder and Sobolev inequalities, it follows from Proposition \ref{proposition 2.10} that
\begin{align}
\label{3.14}
\begin{aligned}
\|\nabla u\|_6 \leq& C\|\rho \theta\|_6+ C\|\rho\dot{u}\|_2+C\|\Delta d\cdot\nabla d\|_2\\
\leq&  C\|\nabla \theta\|_2+ C\|\rho\dot{u}\|_2+C\|\nabla d\|_6\|\nabla^2 d\|_2^{\frac{1}{2}}\|\nabla^2 d\|_6^{\frac{1}{2}}\\
\leq& C\|\nabla \theta\|_2+ C\|\rho\dot{u}\|_2+C\|\nabla^3 d\|_2^{\frac{1}{2}}.
\end{aligned}
\end{align}
Thus, by the Young and Sobolev inequalities, one has
\begin{align}
\label{3.15}
\|\nabla u\|_4^4 \leq C\|\nabla u\|_2\|\nabla u\|_6^3\leq C\|\nabla u\|_2^4+C\|\nabla u\|_6^4
\leq C\big(1+ \|\nabla \theta\|_2^4+ \|\sqrt{\rho}\dot{u}\|_2^4+\|\nabla^3 d\|_2^2\big).
\end{align}
Substituting (\ref{3.15}) into (\ref{3.9}), we obtain
\begin{align}
\label{3.16}
\begin{aligned}
&2\frac{d}{dt}\int \big(\Phi+\eta_1^{\frac{1}{2}}|\nabla d_t|^2 + \eta_1^{\frac{1}{4}}\rho| \dot{u}|^2\big)dx
 + c_v\|\sqrt{\rho}\dot{\theta}\|_2^2+ \eta_1^{\frac{1}{2}}\|(d_{tt}, \Delta d_t)\|_2^2+ \mu\eta_1^{\frac{1}{4}}\|\nabla\dot{u}\|_2^2\\
\leq &
C\big(1+\|\nabla\theta\|_2^2)(1+\|\nabla^3d\|_2^2+\|\nabla\theta\|_2^2\big)
+
C\big(\|\sqrt{\rho}\dot{u}\|_2^4 + \|\nabla d_t\|_2^4+ \|\nabla d_t\|_2^2 + \|\nabla u\|_6^2),
\end{aligned}
\end{align}
where $\Phi$ is given by (\ref{3.7-1}).

Now, we want to show the lower bound of $\Phi$. By the elliptic estimates, it follows from (\ref{1.4}), the H\"{o}lder, Sobolev and Young inequalities that
\begin{align*}
\|\nabla^3 d\|_2
\leq& C\big(\|\nabla d_t\|_2+\|\nabla u\cdot \nabla d\|_2+\|u\cdot\nabla^2d\|_2+\||\nabla d|^2\nabla d\|_2+ \||\nabla d|\nabla^2d\|_2\big)\\
\leq & C\big(\|\nabla d_t\|_2+\|\nabla u\|_4\|\nabla d\|_4+\|u\|_6\|\nabla^2 d\|_2^{\frac{1}{2}}\|\nabla^3 d\|_2^{\frac{1}{2}}+\|\nabla d\|_6^3+
\|\nabla d\|_6\|\nabla^2 d\|_2^{\frac{1}{2}}\|\nabla^3 d|\|_2^{\frac{1}{2}}\big)\\
\leq & \frac{1}{4}\|\nabla^3 d|\|_2+C\big(1+\|\nabla d_t\|_2+\|\nabla u\|_4\big),
\end{align*}
where we have used Proposition \ref{proposition 2.10}.
By (\ref{3.14}) and the Cauchy inequality, it holds that
\begin{align*}
\|\nabla u\|_6 \leq C\big(1+ \|\nabla \theta\|_2+ \|\rho\dot{u}\|_2\big)+ \frac{1}{4}\|\nabla^3 d\|_2.
\end{align*}
Combining the above two inequalities, together with (\ref{3.15}), lead to
\begin{align}
\label{3.17}
\|\nabla u\|_6+\|\nabla^3 d\|_2 \leq
C\big(1+ \|\nabla \theta\|_2+ \|\rho\dot{u}\|_2+\|\nabla d_t\|_2\big).
\end{align}
Thus, from the definition of $\Phi$, (\ref{3.17}), the Young and Sobolev inequalities, one deduces by Proposition \ref{proposition 2.10} that
\begin{align}
\label{3.18}
\begin{aligned}
&2\int \left(\Phi + \eta_1^{\frac{1}{4}}\rho|\dot{u}|^2 +\eta_1^{\frac{1}{2}}|\nabla d_t|^2 \right)dx\\ \geq&
2\kappa \|\nabla\theta\|_2^2 -C\|\theta\|_6\|\nabla u\|_2^{\frac{3}{2}}\|\nabla u\|_6^{\frac{1}{2}}
-C\|\theta\|_6\|\Delta d\|_2^{\frac{3}{2}}\|\Delta d\|_6^{\frac{1}{2}}-C\|\theta\|_6 \|\nabla d \|_{\frac{24}{5}}^4\\
&
+ 2\int\left(\eta_1^{\frac{1}{4}}\rho|\dot{u}|^2 +\eta_1^{\frac{1}{2}}|\nabla d_t|^2\right)dx\\
\geq& \frac{3}{2}\kappa \|\nabla\theta\|_2^2 -C\big(1+\|\nabla u\|_6+ \|\nabla^3d\|_2\big)
+ 2\int\left(\eta_1^{\frac{1}{4}}\rho|\dot{u}|^2 +\eta_1^{\frac{1}{2}}|\nabla d_t|^2\right)dx\\
\geq& \frac{3}{2}\kappa \|\nabla\theta\|_2^2 -C\big(\|\nabla\theta\|_2+\|\sqrt{\rho }\dot{u}\|_2+ \|\nabla d_t\|_2\big)
+ 2\int\left(\eta_1^{\frac{1}{4}}\rho|\dot{u}|^2 +\eta_1^{\frac{1}{2}}|\nabla d_t|^2\right)dx\\
\geq&
\kappa \|\nabla\theta\|_2^2-C(\eta_1) +\int\left(\eta_1^{\frac{1}{4}}\rho|\dot{u}|^2 +\eta_1^{\frac{1}{2}}|\nabla d_t|^2\right)dx.
\end{aligned}
\end{align}

Finally, integrating (\ref{3.16}) over $[0,t]$, and then using (\ref{3.18}) and Gr\"{o}nwall's inequality, the conclusion follows.
\hfill$\Box$

As a straightforward consequence of Lemma \ref{lemma 3.1}, and using (\ref{3.17}),
we have the following corollary:
\begin{corollary}
\label{corollary 3.1}
Assume $\mathcal{N}_0\leq \varepsilon_0$. It holds that
\begin{align*}
\sup_{0\leq t\leq T}(\|\nabla u\|_6+\|\nabla^3 d\|_2+\||\nabla d||\nabla^2 d|\|_2)\leq C_T,
\end{align*}
where $C_T$ depending only on $R$, $c_v$, $\mu$, $\lambda$, $\kappa$, $T$ and the initial data.
\end{corollary}

\begin{lemma}
\label{lemma 3.2}
Assume $\mathcal{N}_0\leq \varepsilon_0$. It holds that
\begin{align*}
\sup_{0\leq t\leq T}(\|\sqrt\rho\dot\theta\|_2^2+\| \nabla^2\theta\|_2^2)+\int_0^T\|\nabla\dot\theta\|_2^2dt\leq C_T,
\end{align*}
where $C_T$ depending only on $R, c_v, \mu, \lambda, \kappa, \Phi_0,$ and $T$.
\end{lemma}

{\it\bfseries Proof.}
Recalling (\ref{3.4}), by the Sobolev inequality, Proposition \ref{proposition 2.10} and Lemma \ref{lemma 3.1}, in order to get this result, it remains to bound the term $\displaystyle\sup_{0\leq t\leq T}\|\sqrt{\rho}\dot{\theta}\|_2^2$.
Applying the operator $\partial_t+{\rm div}(u\cdot)$ to (\ref{1.3}), by tedious computations developed in Appendix A, it follows
\begin{align}
\label{3.19}
\begin{aligned}
&c_v\rho(\dot{\theta}_t+ u\cdot \nabla\dot{\theta})\\
=& \kappa\Delta\dot{\theta}+ \kappa \big({\rm div}u\Delta\theta- \partial_i(\partial_iu\cdot\nabla\theta)-\partial_iu\cdot\nabla\partial_i\theta\big)
+\left(\lambda({\rm div}u)^2+\frac{\mu}{2}|\nabla u+(\nabla u)^t|^2\right){\rm div}u\\
&
 +R\rho\theta\partial_ku_l\partial_lu_k-R\rho\dot{\theta}{\rm div}u
 -R\rho\theta{\rm div}\dot{u}+2\lambda({\rm div}\dot{u}-\partial_ku_l\partial_lu_k){\rm div}u\\&
 +\mu(\partial_iu_j+\partial_ju_i)
 (\partial_i\dot{u}_j+\partial_j\dot{u}_i-\partial_iu_k\partial_ku_j-
 \partial_ju_k\partial_ku_i) \\
 &+\partial_t(|\Delta d+|\nabla d|^2d|^2)+{\rm div}(|\Delta d+|\nabla d|^2d|^2u).
\end{aligned}
\end{align}
Recalling that $|\Delta d+|\nabla d|^2d|^2=|\Delta d|^2-|\nabla d|^4$ and $\Delta d+ |\nabla d|^2d= \Delta d-(d\cdot\Delta d)d$ guaranteed by $|d|=1$, one has $\partial_t|\Delta d+|\nabla d|^2d|^2=2\Delta d\cdot\Delta d_t-4|\nabla d|^2\nabla d:\nabla d_t$. Thanks to this,
multiplying (\ref{3.19}) by $\dot{\theta}$, using integration by parts,
Proposition \ref{proposition 2.10}, Lemma \ref{lemma 3.1}, and  Corollary \ref{corollary 3.1}, we have
\begin{align*}
&\frac{c_v}{2}\frac{d}{dt}\|\sqrt{\rho}\dot{\theta}\|_2^2 + \kappa\|\nabla\dot{\theta}\|_2^2\\
\leq&
C\int |\nabla u|(|\nabla^2\theta||\dot{\theta}|+ |\nabla\theta||\nabla\dot{\theta}|)dx+ \int |\nabla u|^2|\dot{\theta}|(|\nabla u|+\theta)dx\\
&+C\int \rho|\dot{\theta}|^2|\nabla u|dx +C\int \rho\theta|\nabla \dot{u}||\dot{\theta}|dx +C\int |\nabla u||\nabla\dot{u}||\dot{\theta}|dx\\
&+C\int (|\Delta d||\Delta d_t||\dot{\theta}|+ |\nabla d|^3| |\nabla d_t||\dot{\theta}|+ |\Delta d-(d\cdot\Delta d)d|^2|u||\nabla \dot{\theta}| )dx\\
\leq&
C\|\nabla u\|_3\big(\|\nabla^2 \theta\|_{2}\|\dot{\theta}\|_6+\|\nabla\theta\|_6\|\nabla\dot{\theta}\|_2\big)
+C\|\nabla u\|_3^2\|\dot{\theta}\|_6(\|\nabla u \|_6+\|\theta\|_6)\\
&
+C\|\nabla u\|_3\|\rho\dot{\theta}\|_2\|\dot{\theta}\|_6+
C\|\sqrt{\rho}\theta\|_2^{\frac{1}{2}}\|\theta\|_6^{\frac{1}{2}}\|\nabla\dot{u}\|_2
\|\dot{\theta}\|_6 + C\|\nabla u\|_3\|\nabla\dot{u}\|_2\|\dot{\theta}\|_6\\
&+C\|\Delta d\|_3\|\Delta d_t\|_2\|\dot{\theta}\|_6+C\|\nabla d\|_{\frac{9}{2}}^3\|\nabla d_t\|_6\|\dot{\theta}\|_6
+C\|\Delta d\|_6^2\|u\|_6\|\nabla\dot{\theta}\|_2\\
\leq& \frac{\kappa}{2}\|\nabla\dot{\theta}\|_2^2 + C\big(1+ \|\nabla^2\theta\|_2^2+ \|\sqrt{\rho}\dot{\theta}\|_2^2
+ \|\nabla\dot{u}\|_2^2+ \|\Delta d_t\|_2^2\big).
\end{align*}
Thanks to (\ref{1.8}), Lemma \ref{lemma 3.1} and  Corollary \ref{corollary 3.1}, applying Gr\"{o}nwall's inequality, we arrive at
\begin{align*}
\sup_{0\leq t\leq T}\|\sqrt{\rho}\dot{\theta}\|_2^2
+\int^T_0 \| \nabla \dot{\theta}\|_2^2dt\leq C_T,
\end{align*}
which completes the proof.
\hfill$\Box$

\begin{lemma}
\label{lemma 3.3}
Assume $\mathcal{N}_0\leq \varepsilon_0$. It holds that
\begin{align*}
\sup_{0\leq t\leq T}(\|\rho\|_{H^1\cap W^{1, {q}}}+\| \nabla^2 u\|_2)\leq C_T,
\end{align*}
where $C_T$ depending only on $R, c_v, \mu, \lambda, \kappa, \Phi_0,$ and $T$.
\end{lemma}

{\it\bfseries Proof.}
By (\ref{3.111}),  Lemma \ref{lemma 3.1}  and Lemma \ref{lemma 3.2}, it follows
\begin{align}
\label{3.20}
\|\nabla v\|_\infty \leq C_T\log(e+\|\nabla\rho\|_{{q}}), ~~{q}\in (3,6].
\end{align}
Meanwhile, it follows from (\ref{3.13}), Lemma \ref{lemma 3.1} and  Corollary \ref{corollary 3.1}, that
\begin{align}
\label{3.21}
\|\nabla w\|_6 + \|\nabla^2 w\|_2 \leq C\|\sqrt{\rho}\dot{u}\|_2 +C \|\nabla^2d\|_2^{\frac{3}{2}}\|\nabla^3 d\|_2^{\frac{1}{2}}
\leq C_T.
\end{align}
Due to (\ref{3.222}), Proposition \ref{proposition 2.10}, and  Corollary \ref{corollary 3.1}, one deduces by the Sobolev inequality that
\begin{align}
\label{3.22}
\begin{aligned}
\|\nabla^2 w\|_6 \leq C\big(\|\nabla \dot{u}\|_2+ \|\nabla^3d\cdot\nabla d\|_2+ \||\nabla^2 d|^2\|_2\big)
\leq C\big(\|\nabla \dot{u}\|_2+ \|\nabla d\|_{H^2}^2\big)
\leq C\|\nabla \dot{u}\|_2 + C_T.
\end{aligned}
\end{align}
Hence, by the Sobolev inequality, (\ref{3.21}) and (\ref{3.22}) give us
\begin{align*}
\|\nabla w\|_\infty \leq C\|\nabla \dot{u}\|_2 + C_T,
\end{align*}
which combined with (\ref{3.20}) implies
\begin{align}
\label{3.23}
\|\nabla u\|_\infty \leq C_T\log(e+\|\nabla\rho\|_{{q}}) + C\|\nabla \dot{u}\|_2, ~~{q}\in (3,6].
\end{align}
Applying the elliptic estimates to (\ref{1.2}), one has for $2\leq p\leq {q}$
\begin{align}
\label{3.24}
\begin{aligned}
\|\nabla^2 u\|_p \leq& C (\|\rho\dot{u}\|_p+\|\Delta d\cdot \nabla d\|_p+\|\nabla P\|_p)\\
\leq& C(\|\rho\dot{u}\|_p +\|\Delta d\|_p\|\nabla d\|_\infty+ \|\rho\nabla\theta\|_p+\|\nabla\rho\theta\|_p)\\
\leq& C(\|\rho\dot{u}\|_p +\|\nabla d\|_{H^2}^2+ \|\nabla\theta\|_{H^1}+\|\nabla\rho\|_p\|\theta\|_\infty)\\
\leq& C(1+ \|\rho\dot{u}\|_p+\|\nabla \rho\|_p)\\
\leq& C_T(1 + \|\nabla\dot{u}\|_2+\|\nabla \rho\|_p),
\end{aligned}
\end{align}
where Proposition \ref{proposition 2.10} and Lemma \ref{lemma 3.2} were used.
On the other hand, some straightforward calculations show that, for $2\leq p\leq {q}$
\begin{align}
\label{3.25}
\frac{d}{dt}\|\nabla \rho\|_p \leq C(1+\|\nabla u\|_\infty)\|\nabla\rho\|_p+ C\|\nabla^2 u\|_p,
\end{align}
which together with (\ref{3.23}) and (\ref{3.24}) yields
\begin{align*}
\frac{d}{dt}\|\nabla \rho\|_p \leq C_T\big(1+ \log (e+\|\nabla\rho\|_{{q}})+ \|\nabla\dot{u}\|_2\big)\|\nabla\rho\|_p
 + C_T(1 + \|\nabla\dot{u}\|_2+\|\nabla \rho\|_p).
\end{align*}
Set
$$f(t)= e+\|\nabla\rho\|_{{q}} \text{~~ and ~~} g(t)=1+\|\nabla\dot{u}\|_2,$$
then
\begin{align*}
\frac{d}{dt}f(t) \leq C_Tg(t)f(t)\log f(t).
\end{align*}
By solving the above ordinary differential inequality and using Lemma \ref{lemma 3.1}, one gets
\begin{align}
\label{3.26}
\sup_{0\leq t\leq T}\|\nabla\rho\|_{{q}}\leq C_T.
\end{align}
Combing (\ref{3.26}) with (\ref{3.23}) yields
\begin{align}
\label{3.27}
\int ^T_0 \|\nabla u\|_\infty^2dt \leq C_T.
\end{align}
Choosing $p=2$ in (\ref{3.25}), it follows from Lemma \ref{lemma 3.1}, (\ref{3.24}), (\ref{3.27}) and Gr\"{o}nwall's inequality that
\begin{align*}
\sup_{0\leq t\leq T}\|\nabla\rho\|_{2}\leq C_T.
\end{align*}
This together with (\ref{3.24}) and Lemma \ref{lemma 3.1} implies
\begin{align*}
\sup_{0\leq t\leq T}\|\nabla^2 u\|_{2}\leq C(\|\rho\dot{u}\|_2+\|\nabla \rho\|_2) + C_T\leq C_T.
\end{align*}
Thus this lemma is proved.
\hfill$\Box$

As a direct consequence of Lemma \ref{lemma 3.1}, Lemma \ref{lemma 3.2}, Lemma \ref{lemma 3.3}, and Corollary \ref{corollary 3.1}, we have the following Corollary.

\begin{corollary}
  \label{corollary 3.2}
  Assume that $\mathcal N_0\leq\varepsilon_0$. Then, it holds that
  \begin{eqnarray*}
    &\displaystyle\sup_{0\leq t\leq T}(\|\nabla\theta\|_{H^1}^2+\|(\nabla^2u, \sqrt\rho\dot u, \sqrt\rho\dot\theta, \nabla^3d, \nabla d_t)\|_2^2+\|\rho\|_{H^1\cap W^{1,q}})\leq C_T, \\
    &\displaystyle\int_0^T\|(\nabla\dot u, d_{tt}, \nabla^2d_t, \nabla\dot\theta)\|_2^2dt\leq C_T,
  \end{eqnarray*}
  for a positive constant $C_T$ depending only on $R, c_v, \mu, \lambda, \kappa, \Phi_0,$ and $T$.
\end{corollary}

\section{Proof of Theorem \ref{theorem 1.1}}

{\it\bfseries Proof of Theorem \ref{theorem 1.1}.}
Let $(\rho, u, \theta, d)$ be the unique local solution guaranteed by Lemma \ref{lemma 2.1}.
By applying the local well-posedness, i.e.\,Lemma \ref{lemma 2.1}, inductively, one can extend the $(\rho, u, \theta, d)$
uniquely
to the maximal time $T_{\text{max}}$ of existence. We claim that $T_\text{max}=\infty$ and thus the conclusion holds.
Assume by contradiction that $T_\text{max}<\infty$. Let $\varepsilon_0$ be as in Lemma \ref{lemma 2.9} and assume that
$\mathcal{N}_0\leq \varepsilon_0$. Then, it follows from Proposition \ref{proposition 2.10} and Corollary \ref{corollary 3.2} that
$$
\sup_{0\leq t\leq T}\left(\|\rho\|_{W^{1,q}\cap H^1}(t)+\|(u,\theta)\|_{D_0^1\cap D^2}(t)+\|\nabla d\|_{H^2}(t)+\|(\sqrt{\rho}\theta, \sqrt\rho\dot u, \sqrt\rho\dot\theta)\|_2(t)\right)\leq C_0,
$$
for any $T\in(0,T_\text{max})$, where $C_0$ is a positive constant independent of $T\in(0,T_\text{max})$. Take
$\delta>0$ and choose $T_\text{max}-\delta$ as the initial time. Thanks to the above estimate, one can check that all the conditions in Lemma \ref{lemma 2.1} hold and, in particular, the compatibility conditions hold.
Then, by Lemma \ref{lemma 2.1}, there is a positive time $\mathcal T_0$, depending only on
$R, c_v, \mu, \lambda, \kappa,$
and $C_0$ (thus independent of $\delta$), such that the solution $(\rho, u, \theta, d)$ can be extended to the time
$T_\text{max}-\delta+\mathcal T_0$. By choosing $\delta$ sufficiently small, it holds that
$T_\text{max}-\delta+\mathcal T_0>T_\text{max}$, which contradicts to the definition of $T_\text{max}$.
This contradiction implies that $T_\text{max}=\infty$. Therefore, we get a global strong solution and the conclusion
follows.
\hfill$\Box$

\section*{Acknowledgement}
{\small J. Li was supported in part by the National
Natural Science Foundation of China (11971009 and 11871005), by the Key Project of National Natural Science Foundation of China (12131010), and by the Guangdong Basic and Applied Basic Research Foundation (2019A1515011621,
2020B1515310005, 2020B1515310002, and 2021A1515010247).
Q. Tao was supported in part by the National Science Foundation
of China (11971320), and by the Guangdong Basic and Applied Basic Research Foundation (2020A1515010530).}

\section*{Appendix A}
Applying the operator $\partial_t+{\rm div}(u\cdot)$ to (\ref{1.3}), from the definition of material derivative, one has
\begin{align}
\label{4.1}
\begin{aligned}
&c_v [\partial_t(\rho\dot{\theta})+{\rm div}(u \rho\dot{\theta})] +\partial_t(P{\rm div}u)+ {\rm div}(u P{\rm div}u)-\kappa[\partial_t\Delta\theta+ {\rm div}(u\Delta\theta)]\\
=& \frac{\mu}{2}\left[\partial_t\left(|\partial_ju_i+\partial_iu_j|^2\right)+ {\rm div}(u |\partial_ju_i+\partial_iu_j|^2)\right]+
\lambda[\partial_t({\rm div}u)^2+ {\rm div}(u({\rm div}u)^2)]\\
&+\partial_t(|\Delta d+|\nabla d|^2d|^2)+{\rm div}(|\Delta d+|\nabla d|^2d|^2u).
\end{aligned}
\end{align}
By (\ref{1.1}) and some straightforward calculations, it follows that
\begin{align}
&c_v [\partial_t(\rho\dot{\theta})+{\rm div}(u \rho\dot{\theta})]= c_v\rho(\dot{\theta}_t+ u\cdot \nabla\dot{\theta}),\label{4.2}\\[3mm]
&\partial_t(P{\rm div}u)+ {\rm div}(u P{\rm div}u)
=R \big[\partial_t\rho\theta{\rm div}u+ \rho\partial_t\theta{\rm div}u+ \rho\theta{\rm div}u_t+\rho\theta({\rm div}u)^2\nonumber\\
&\qquad\qquad\qquad\qquad\qquad\qquad+ u\cdot\nabla\rho\theta{\rm div}u+ \rho u\cdot\nabla\theta{\rm div}u
+\rho\theta u\cdot\nabla( {\rm div}u)\big]\nonumber\\
&\qquad\qquad\qquad\qquad\qquad\quad=R \big[ \rho\dot{\theta}{\rm div}u+ \rho\theta {\rm div}u_t+ \rho\theta u\cdot\nabla( {\rm div}u)\big]\label{4.3}\\
&\qquad\qquad\qquad\qquad\qquad\quad= R \big[ \rho\dot{\theta}{\rm div}u+ \rho\theta {\rm div}\dot{u} -\rho\theta{\rm div}(u\cdot\nabla u) + \rho\theta u\cdot\nabla( {\rm div}u)\big]\nonumber\\
&\qquad\qquad\qquad\qquad\qquad\quad= R (\rho\dot{\theta}{\rm div}u+ \rho\theta {\rm div}\dot{u} )-R\rho\theta\partial_ku_l\partial_lu_k,\nonumber\\[3mm]
&-\kappa[\partial_t\Delta\theta+ {\rm div}(u\Delta\theta)]= -\kappa[\Delta \dot{\theta} - \Delta(u\cdot\nabla\theta)+ {\rm div}(u\Delta\theta)]\nonumber\\
&\qquad\qquad\qquad\qquad\quad~~= -\kappa[\Delta \dot{\theta} - \Delta(u\cdot\nabla\theta)+ {\rm div}u\Delta\theta+ u\cdot\nabla(\Delta\theta)]\label{4.4}\\
&\qquad\qquad\qquad\qquad\quad~~= -\kappa\Delta \dot{\theta} -\kappa\big[{\rm div}u\Delta\theta - \partial_i(\partial_iu\cdot\nabla\theta)-\partial_iu\cdot\nabla\partial_i\theta\big].\nonumber
\end{align}
Similarly, for the terms on the right hand side of (\ref{4.1}), we also have
\begin{align}
\label{4.5}
\begin{aligned}
&\frac{\mu}{2}\left[\partial_t\left(|\partial_ju_i+\partial_iu_j|^2\right)+ {\rm div}(u |\partial_ju_i+\partial_iu_j|^2\right]\\
=&
\mu(\partial_ju_i+\partial_iu_j)(\partial_t\partial_ju_i+\partial_t\partial_iu_j)+\frac{\mu}{2}|\nabla u+(\nabla u)^t|^2{\rm div}u
+\mu u\cdot(\partial_ju_i+\partial_iu_j)\nabla(\partial_ju_i+\partial_iu_j)
\\
=&\frac{\mu}{2}|\nabla u+(\nabla u)^t|^2{\rm div}u+ \mu(\partial_iu_j+\partial_ju_i)
 (\partial_i\dot{u}_j+\partial_j\dot{u}_i-\partial_iu_k\partial_ku_j-
 \partial_ju_k\partial_ku_i),
\end{aligned}
\end{align}
and
\begin{align}
\label{4.6}
\begin{aligned}
&\lambda[\partial_t({\rm div}u)^2+ {\rm div}(u({\rm div}u)^2)]\\
=&2\lambda(\partial_t{\rm div}u){\rm div}u+2\lambda(u\cdot\nabla{\rm div}u){\rm div}u+ \lambda({\rm div}u)^3\\
=&2\lambda {\rm div}\dot{u} {\rm div}u-2\lambda{\rm div}(u\cdot\nabla u){\rm div}u+2\lambda(u\cdot\nabla{\rm div}u){\rm div}u+ \lambda({\rm div}u)^3\\
=&\lambda({\rm div}u)^3+2\lambda({\rm div}\dot{u}-\partial_ku_l\partial_lu_k){\rm div}u.
\end{aligned}
\end{align}
Thus, (\ref{3.19}) follows by combing (\ref{4.1})-(\ref{4.6}) together.


\vskip8mm
\begin{thebibliography}{10}

\bibitem{CK2006} Y. Cho, H. Kim, Existence results for viscous polytropic fluids with vacuum, J. Differential Equations 228 (2006) 377--411.

\bibitem{DeAnnaLiu}
F. De Anna, C. Liu, Non-isothermal general Ericksen-Leslie system: derivation, analysis and thermodynamic consistency,
Arch. Ration. Mech. Anal., 231 (2019), no. 2, 637--717.

\bibitem{DWW2011} S. Ding, C. Wang, H. Wen, Weak solution to compressible hydrodynamic flow of liquid
crystals in dimension one, Discrete Contin. Dyn. Syst. Ser. B, 15 (2011) 357--371.

\bibitem{DLWW2012} S. Ding, J.  Lin, C. Wang, H. Wen, Compressible hydrodynamic flow of liquid crystals
in 1-D, Discrete Contin. Dyn. Syst., 32 (2012) 539--563.


\bibitem{E1962} J. Ericksen, Hydrostatic theory of liquid crystal, Arch. Rat. Mech. Anal., 9 (1962) 371--378.




\bibitem{FLN2018} J. Fan, F. Li, G. Nakamura, Local well-posedness for a compressible non-isothermal model for nematic liquid crystals,
J. Math. Phys. 59 (2018) 031503



\bibitem{FFRS2012} E. Feireisl, M. Fr\'{e}mond, E. Rocca, G. Schimperna, A new approach to non-isothermal models for nematic liquid crystals, Arch. Ration. Mech. Anal. 205 (2012) 651--672.

\bibitem{FRS2011} E. Feireisl, E. Rocca, G. Schimperna, On a non-isothermal model for nematic liquid crystals, Nonlinearity 24 (2011) 243--257.

\bibitem{GTY2016} J. Gao, Q. Tao, Z. Yao, Long-time behavior of solution for the compressible nematic liquid crystal flows in $\mathbb{R}^3$, J. Differential Equations 261 (2016) 2334--2383.

\bibitem{NONUNIQUENESS}
H. Gong, T. Huang, J. Li, Nonuniqueness of nematic liquid crystal flows in dimension three,
J. Differential Equations, 263 (2017), no. 12, 8630--8648.

\bibitem{GXX2017} B. Guo, X. Xi, B. Xie, Global well-posedness and decay of smooth solutions to the non-isothermal model for compressible nematic liquid crystals, J. Differential Equations 262 (2017) 1413--1460.

\bibitem{HONGWEAK}
M.-C. Hong, Global existence of solutions of the simplified Ericksen-Leslie system in dimension two, Calc. Var. Partial Differential Equations, 40 (2011), no. 1-2, 15--36.

\bibitem{HONGXINWEAK}
M.-C. Hong, Z. Xin, Global existence of solutions of the liquid crystal flow for the Oseen-Frank model in $\mathbb R^2$,
Adv. Math., 231 (2012), no. 3-4, 1364--1400.

\bibitem{HLZ2014} M. Hong, J. Li, Z. Xin, Blow-up criteria of strong solutions to the Ericksen-Leslie
system in $\mathbb{R}^3$, Comm. Partial Differential Equations, 39 (2014) 1284--1328.

\bibitem{HW2013} X. Hu, H. Wu, Global solution to the three-dimensional compressible flow of liquid crystals,
SIAM J. Math. Anal., 45 (2013) 2678--2699.

\bibitem{BLOWUP}
T. Huang, F. Lin, C. Liu, C. Wang, Finite time singularity of the nematic liquid crystal flow in dimension three,
Arch. Ration. Mech. Anal., 221 (2016), no. 3, 1223--1254.

\bibitem{HW2012} T. Huang, C. Wang, Blow up criterion liquid crystal flows, Comm. Partial Differential Equations 37 (2012)
875--884.

\bibitem{HWW2012} T. Huang, C. Wang, H. Wen, Strong solutions of the compressible nematic liquid crystal flow,
J. Differential Equations 252 (2012) 2222--2256.

\bibitem{HWW2012-2}  T. Huang, C. Wang, H. Wen, Blow up criterion for compressible nematic liquid crystal
flows in dimension three, Arch. Rat. Mech. Anal., 204 (2012) 285--311.

\bibitem{HL2013} X. Huang, J. Li, Serrin-type blowup criterion for viscous,
compressible, and heat conducting Navier-Stokes
and magnetohydrodynamic flows, Commun. Math. Phys. 324 (2013) 147--171.

\bibitem{HL2018} X. Huang, J. Li, Global classical and weak solutions to the three-dimensional full compressible
Navier-Stokes system with vacuum and large oscillations, Arch. Rational Mech. Anal. 227 (2018) 995--1059.


\bibitem{JJW2013} F. Jiang, S. Jiang, D. Wang, On multi-dimensional compressible flows of nematic liquid crystals with large initial energy in a bounded domain, J. Funct. Anal. 265 (2013) 3369--3397.

\bibitem{JJW2014} F. Jiang, S. Jiang, D. Wang, Global weak solutions to the equations of compressible flow of nematic liquid crystals in two dimensions, Arch. Ration. Mech. Anal. 214 (2014) 403--451.



\bibitem{LLZ2014} Z. Lei, D. Li, X. Zhang, Remarks of global wellposedness of liquid crystal flows and heat flows of harmonic maps in two dimensions, Proc. Am. Math. Soc. 142 (2014) 3801--3810.

\bibitem{L1968} F. Leslie, Some constitutive equations for liquid crystals, Arch. Rat. Mech. Anal., 28 (1968) 265--283.


\bibitem{L2020} J. Li, Global small solutions of heat conductive compressible Navier-Stokes equations with
vacuum: smallness on scaling invariant quantity, Arch. Rational Mech. Anal. 237 (2020) 899--919.

\bibitem{LITITIXIN}
J. Li, E. S. Titi, Z. Xin, On the uniqueness of weak solutions to the Ericksen-Leslie liquid crystal model in $\mathbb R^2$, Math. Models Methods Appl. Sci., 26 (2016), no. 4, 803--822.

\bibitem{LIXINNONISO}
J. Li, Z. Xin, Global existence of weak solutions to the non-isothermal nematic liquid crystals in 2D, Acta Math. Sci. Ser. B (Engl. Ed.), 36 (2016), no. 4, 973--1014.

\bibitem{LXZ2018} J. Li, Z. Xu, J. Zhang, Global existence of classical solutions with large oscillations and vacuum to the three-dimensional compressible nematic liquid crystal flows, J. Math. Fluid Mech. 20 (2018) 2105--2145.


\bibitem{L1989} F. Lin, Nonlinear theory of defects in nematic liquid crystals: phase transition and flow phenomena, Commun. Pure Appl. Math. 42 (1989) 789--814.

\bibitem{LINLINWANG}
F. Lin, J. Lin, C. Wang, Liquid crystal flows in two dimensions, Arch. Ration. Mech. Anal., 197 (2010), no. 1, 297--336.

\bibitem{LL1995}
F. Lin, C. Liu, Nonparabolic dissipative systems modeling the flow of liquid crystals, Commun. Pure Appl. Math.
48 (1995) 501--537.

\bibitem{LL1996}
F. Lin, C. Liu, Partial regularity of the dynamic system modeling the flow of liquid cyrstals, Discrete Contin. Dyn.
Syst. 2 (1996) 1--22.

\bibitem{LINWANGUNIQUENESS}
F. Lin, C. Wang, On the uniqueness of heat flow of harmonic maps and hydrodynamic flow of nematic liquid crystals,
Chin. Ann. Math. Ser. B, 31 (2010), no. 6, 921--938.


\bibitem{LINWANGWEAK3D}
F. Lin, C. Wang, Global existence of weak solutions of the nematic liquid crystal flow in dimension three, Commun. Pure Appl. Math., 69 (2016), no. 8, 1532--1571.

\bibitem{LLW2015}
J. Lin, B. Lai, C. Wang, Global finite energy weak solutions to the compressible nematic liquid crystal flow
in dimension three, SIAM J. Math. Anal. 47 (2015) 2952--2983.

\bibitem{LZ2021} Y. Liu, X. Zhong, Global well-posedness to the 3D Cauchy problem of compressible non-isothermal nematic liquid crystal flows with vacuum,
 Nonlinear Anal. Real World Appl.  58  (2021), 103219, 24 pp.



\bibitem{GTY2015} Q. Tao, J. Gao, Z. Yao, Global strong solutions of the compressible nematic liquid crystal flow with the cylinder symmetry, Commun. Math. Sci.  13  (2015) 2065--2096.



\bibitem{W2016} T. Wang, Global existence and large time behavior of strong solutions to the 2-D compressible nematic liquid crystal flows with vacuum, J. Math. Fluid Mech. 18 (2016) 539--569.

\bibitem{WANGWANGZHANGUNIQUENESS}
M. Wang, W. Wang, Z. Zhang, On the uniqueness of weak solution for the 2-D Ericksen-Leslie system,
Discrete Contin. Dyn. Syst. Ser. B, 21 (2016), no. 3, 919--941.


\bibitem{XZ2012} X. Xu, J. Zhang, A blow-up criterion for 3D compressible magnetohydrodynamic equations with
vaccum, Math. Models Method. Appl. Sci. 22 (2012), 1150010.

\bibitem{ZX2019} X. Zhong, Singularity formation to the two-dimensional compressible non-isothermal nematic liquid crystal flows in a bounded domain, J. Differential Equations 267(2019) 3797--3826.

\end{thebibliography}
\end{document}